\documentclass[a4paper,12pt]{article}

\usepackage[margin=3cm,footskip=1cm]{geometry}

\usepackage[T1]{fontenc}

\usepackage{amsmath,amssymb,amsthm,bm,mathtools,xspace,booktabs,url}
\usepackage{graphicx,etoolbox}

\usepackage[shortlabels]{enumitem}
\setlist{
  listparindent=\parindent,
  parsep=0pt,
}

\usepackage[round]{natbib}

\AtBeginEnvironment{thebibliography}{%
  \small
  \setlength\itemsep{0.75em plus 0.5em minus 0.75em}%
}

\usepackage{caption}
\usepackage[raggedright]{titlesec}

\numberwithin{equation}{section}

\usepackage{authblk}

\makeatletter
\newcommand\CorrespondingAuthor[1]{%
  \begingroup%
  \def\@makefnmark{}%
  \footnotetext{Corresponding author: #1}%
  \endgroup%
}
\makeatother

\makeatletter
\renewenvironment{abstract}{%
  \small%
  \providecommand\keywords{%
    \par\medskip\noindent\textit{Keywords:}\xspace}%
  \begin{center}%
    \bfseries \abstractname\vspace{-.5em}\vspace{\z@}%
  \end{center}%
  \quote%
}{\endquote}
\makeatother

\usepackage[above,below,section]{placeins}

\usepackage[load-configurations=abbreviations]{siunitx}
\usepackage{lipsum}

\newtheoremstyle{example}
{\topsep} {\topsep}%
{\upshape}
{}
{\itshape}
{.}
{1em}
{\thmname{#1}\thmnumber{ #2 }\thmnote{#3}}

\theoremstyle{example}
\newtheorem{example}{Example}
\newtheorem*{myexample}{}

\newcommand\psup[1]{^{\smash{(}#1\smash{)}}}

\begin{document}

\title{Structured space-sphere point processes and $K$-functions}

\author{Jesper M\o ller}

\author{Heidi S. Christensen}

\author{Francisco Cuevas-Pacheco}

\author{Andreas D. Christoffersen}

\affil{Department of Mathematical Sciences, Aalborg University}

%
%

\date{}


\maketitle

\begin{abstract}

  This paper concerns space-sphere point processes, that is, point
  processes on the product space of $\mathbb R^d$ (the $d$-dimensional
  Euclidean space) and $\mathbb S^k$ (the $k$-dimen\-sional
  sphere). We consider specific classes of models for space-sphere
  point processes, which are adaptations of existing models for either
  spherical or spatial point processes.  For model checking or
  fitting, we present the space-sphere $K$-function which is a natural
  extension of the inhomogeneous $K$-function for point processes on
  $\mathbb R^d$ to the case of space-sphere point processes.  Under
  the assumption that the intensity and pair correlation function both
  have a certain separable structure, the space-sphere $K$-function is
  shown to be proportional to the product of the inhomogeneous spatial
  and spherical $K$-functions. For the presented space-sphere point
  process models, we discuss cases where such a separable structure
  can be obtained.  The usefulness of the space-sphere $K$-function is
  illustrated for real and simulated datasets with varying dimensions
  $d$ and $k$.


\end{abstract}

\section{Introduction}\label{intro}

Occasionally point processes arise on more complicated spaces than the usual space $\mathbb{R}^d$, the $d$-dimensional Euclidean space, as for spatio-temporal point processes, spherical point processes or point processes on networks \citep[see][and the references therein for details on such point processes]{PD2016,LBMN2016,MR2016,BNMR2017}. In this paper we consider \textit{space-sphere point processes} that live on the product space $S = \mathbb{R}^d \times \mathbb{S}^k$, where $\mathbb S^k=\{u\in\mathbb R^{k+1}:\|u\|_{k+1}=1\}$ is the $k$-dimensional unit sphere, $\|\cdot\|_k$ denotes the usual distance in $\mathbb R^k$, and  $d,k\in\{1,2,\ldots\}$. For each point $(y, u) \in S$ belonging to a given space-sphere point process, we call $y$ its spatial component and $u$ its spherical component. Assuming local finiteness of a space-sphere point process, the spatial components constitute a locally finite point process in $\mathbb{R}^d$, but the spherical components do not necessarily form a finite point process on  $\mathbb{S}^k$. However, in practice the spatial components are only considered within a bounded window $W \subset \mathbb{R}^d$, and the associated spherical components do constitute a finite point process. 

\begin{figure}[htp]	
	\centering
	\includegraphics[width=.45\textwidth, trim = {0cm 2cm 0cm 2cm}, clip]{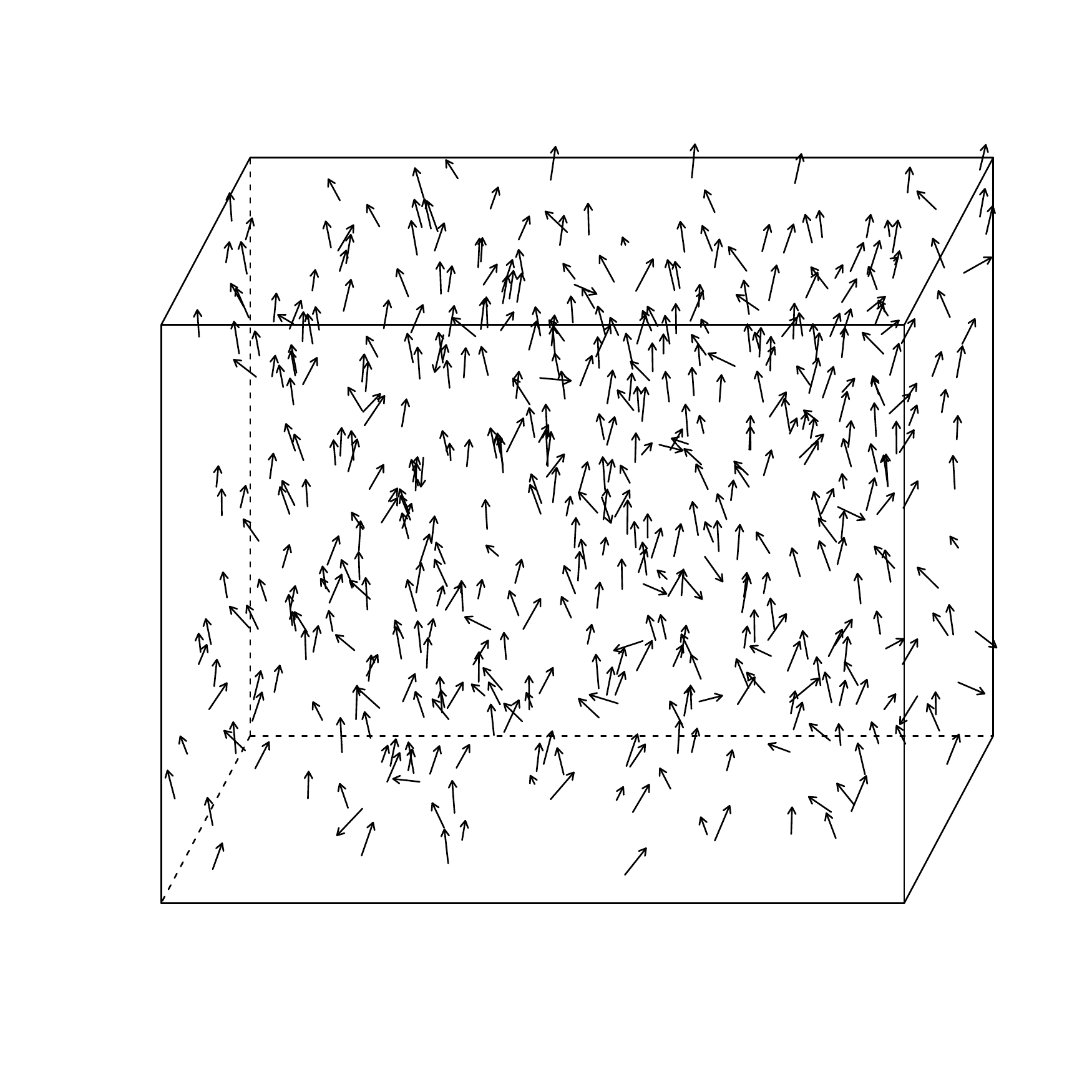}
	\caption{Location and orientation of pyramidal neurons in a small section of a human brain. For details, see Section~\ref{s:neuron}.}
	\label{fig:brain_cells}
\end{figure}

One example is the data shown in Figure~\ref{fig:brain_cells} that
consists of the location and orientation of a number of pyramidal
neurons found in a small area of a healthy human's primary motor
cortex. More precisely, the locations are three-dimensional
coordinates each describing the placement of a pyramidal neuron's
nucleolus, and the orientations are unit vectors pointing from a
neuron's nucleolus toward its apical dendrite. These data can be
considered as a realisation of a space-sphere point process with
dimensions $d = 3$ and $k = 2$, where the spatial components describe
the nucleolus locations and the spherical components are the
orientations.  How neurons (of which around \SIrange{75}{80}{\%} are
pyramidal neurons) are arranged have been widely discussed in the
literature. Specifically, it is hypothesised that neurons are arranged
in columns perpendicular to the pial surface of the brain. This
hypothesis, referred to as the minicolumn hypothesis, have been
studied for more than half a century \citep[see e.g.][]{Lorente1938,
  Mountcastle:1978, Buxhoeveden2002}, and it is believed that
deviation from such a columnar structure is linked with neurological
diseases such as Alzheimers and schizophrenia.

\begin{figure}[htp]
	\centering
	\includegraphics[width = 0.4\textwidth, angle = -90, trim = {2cm 5cm 2cm 5cm}, clip]{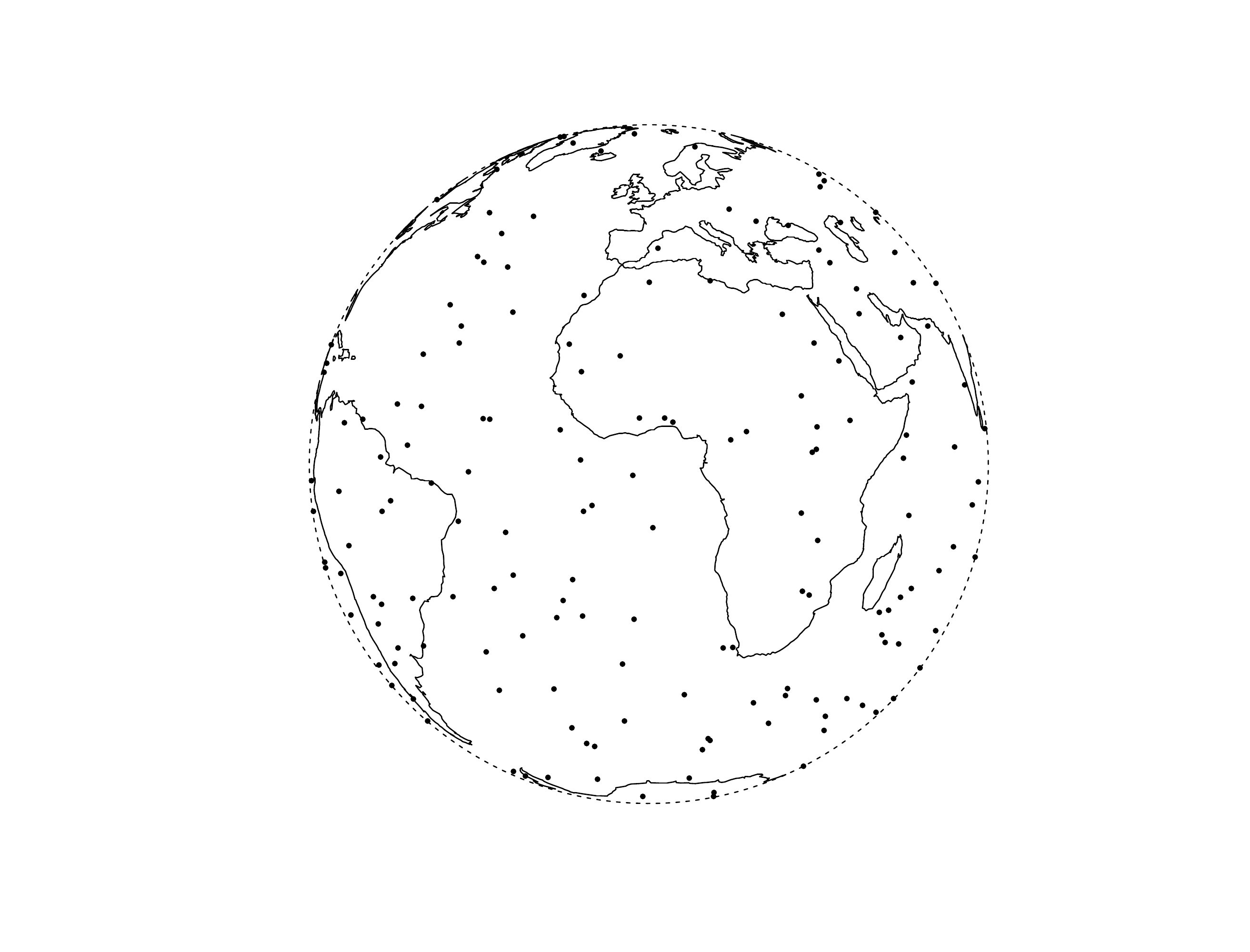}
	\includegraphics[width = 0.4\textwidth, angle = -90, trim = {2cm 5cm 2cm 5cm}, clip]{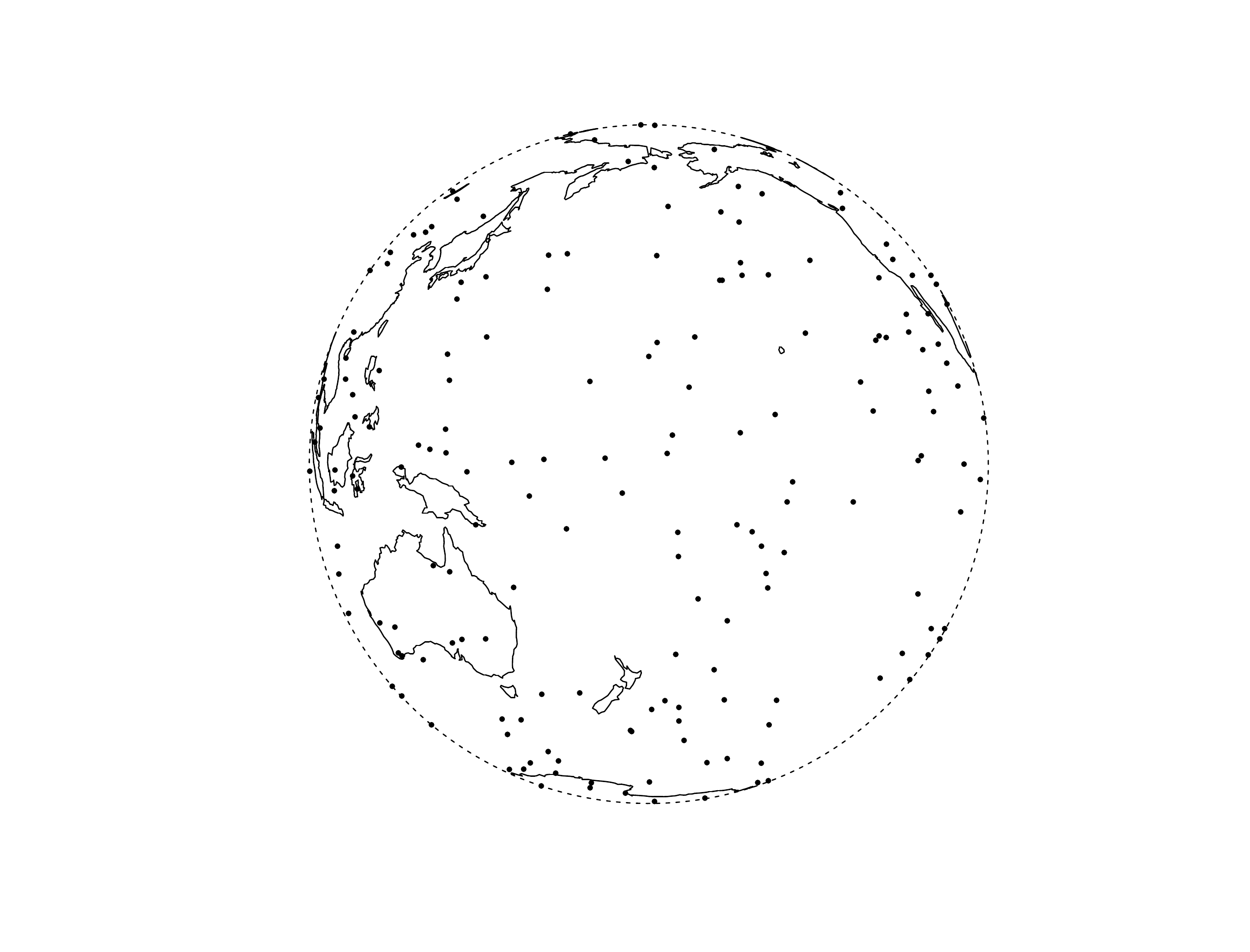}
	\includegraphics[width = 0.1\textwidth, angle = -90]{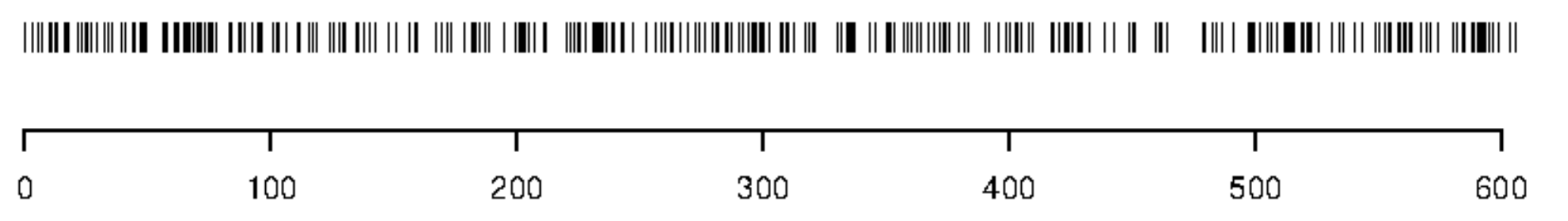}
	\caption{Top: orthographic projection of the fireball locations. Bottom: time of fireball events measured in weeks. For details, see Section~\ref{sec:fireball}.}
	\label{fig:MapFireballsSphere}
\end{figure}

Another example is the time and geographic location of fireballs, which are bright meteors reaching a visual magnitude of $-3$ or brighter. They are continually recorded by U.S.\ Government sensors and made available at \url{http://neo.jpl.nasa.gov/fireballs/}. We can consider fireball events as a space-sphere point process with dimensions $d = 1$ and $k = 2$, where the time and locations are the spatial and spherical components, respectively. 
Figure~\ref{fig:MapFireballsSphere} shows the location of fireballs on the globe (identified with the unit sphere) observed over a time period of about 606 weeks. 

The paper is organised as follows. In Section~\ref{sec:preliminary}, we define concepts related to space-sphere point processes and give some natural examples of such processes. In Section~\ref{s:SOIRS}, we define the space-sphere $K$-function, a functional summary statistic which is analogue to the space-time $K$-function when $d=2$ and $\mathbb{S}^k$ is replaced by the time axis \citep{Diggle1995, Gabriel2009, Moeller2012}.
The space-sphere $K$-function is defined in terms of the pair correlation function which is assumed to have a certain stationary form. In the case where both the intensity and pair correlation function have a specific separable structure discussed in Section~ \ref{s:separability}, the space-sphere $K$-function is shown to be proportional to the product of the spatial $K$-function \citep{BMW2000} and the spherical $K$-function \citep{LBMN2016, MR2016}. Further, an unbiased estimate is given in Section~\ref{s:est}. In Section~\ref{s:appl}, the usefulness of the space-sphere $K$-function is illustrated for the fireball and neuron data as well as for simulated data, and it is e.g.\ seen how the $K$-function may be used to test for independence between the spatial and spherical components. 

\section{Preliminaries}\label{sec:preliminary}
\subsection{Setting}\label{sec:setting}
Throughout this paper we consider the following setting. 

Equip $\mathbb R^d$ with the Lebesgue measure $|A|=\int_A \, \mathrm dy$ and $\mathbb S^k$ with Lebesgue/surface measure $\nu(B)$, where $A\subseteq\mathbb R^d$ and $B\subseteq\mathbb S^k$ are Borel sets. 
Thus, the product space 
$S=\mathbb R^d\times\mathbb S^k$ is equipped with Lebesgue measure $\mu$ given by $\mu(A\times B)=|A|\nu(B)$.

Let $X$ be a simple locally finite point process on $S$, that is, we can view $X$ as
a random subset of $S$ such that the restriction $X_B=X\cap B$ of $X$ to any bounded set $B\subset S$ is finite. 
We call $X$ a \textit{space-sphere point process}, and assume that it has intensity function $\rho$ with respect to $\mu$ and pair correlation function $g$ with respect to the product measure $\mu\otimes\mu$. That is, 
for any Borel function $h: S\mapsto [0, \infty)$, 
\begin{align}\label{eq:defn_intensity}
\mathrm{E} \Bigl\{\sum_{(y_i, u_i) \in X} h(y_i, u_i)\Bigr\}
= \int h(y, u)\rho(y, u)\,\mathrm d\mu(y, u),
\end{align} 
provided this integral is finite. We say that $X$ is \textit{(first order) homogeneous} if $\rho$ is a constant function. Furthermore, for any Borel function $k: S\times S \mapsto [0, \infty)$, 
\begin{equation}\label{eq:defn_2nd_order_intensity}
\begin{aligned}
\mathrm{E}& \Bigl\{\sum^{\neq}_{(y_i, u_i), (y_j, u_j) \in X} k(y_i, u_i, y_j, u_j)\Bigr\} \\
&= \iint  k(y_1, u_1, y_2, u_2)\rho(y_1, u_1)\rho(y_2, u_2)g(y_1, u_1, y_2, u_2) 
\,\mathrm d\mu(y_1, u_1)\,\mathrm d\mu(y_2, u_2),
\end{aligned}
\end{equation}
provided this double integral is finite. Here, we set
$g(y_1, u_1, y_2, u_2)=0$ if $\rho(y_1, u_1)\cdot\rho(y_2, u_2)=0$, and
$\sum^{\neq}_{(y_i, u_i), (y_j, u_j) \in X}$ means that we sum over
pairs of distinct points $(y_i, u_i), (y_j, u_j) \in X$.

The functions $\rho$ and $g$ are unique except for null sets with respect to $\mu$ and $\mu\otimes\mu$, respectively. For ease of presentation, we ignore null sets in the following. Note that $g(y_1,u_1,y_2,u_2)=g(y_2,u_2,y_1,u_1)$ is symmetric on $S\times S$. We say that $X$ is {\it stationary in space} if its distribution is invariant under translations of its spatial components; this implies that $\rho(y,u)$ depends only on $u$, and $g(y_1,u_1,y_2,u_2)$ depends only on $(y_1,y_2)$ through the difference $y_1-y_2$. If the distribution of $X$ is invariant under rotations (about the origin in $\mathbb R^d$) of its spatial components, we say that $X$ is {\it isotropic in space}.
Stationarity and isotropy in space imply that $g(y_1,u_1,y_2,u_2)$ depends only on $(y_1,y_2)$ through the distance $\|y_1-y_2\|_d$. 
We say that $X$ is {\it isotropic on the sphere} if its distribution is invariant under rotations (on $\mathbb S^k$) of its spherical components;
this implies that $g(y_1,u_1,y_2,u_2)$ depends only on $(u_1,u_2)$ through the geodesic (great circle/shortest path) distance $d(u_1,u_2)$ on $\mathbb S^k$. If $X$ is stationary in space and isotropic on the sphere, then $\rho$ is constant and 
\begin{equation}\label{e:gSOIRS}
g(y_1, u_1, y_2, u_2) = g_0\{y_1 - y_2, d(u_1, u_2)\}, \qquad y_1,y_2\in\mathbb R^d,\ u_1,u_2\in\mathbb S^k,
\end{equation} 
depends only on $(y_1, y_2)$ through $y_1-y_2$ and on $(u_1, u_2)$ through $d(u_1, u_2)$ (this property is studied further in Section~\ref{s:SOIRS}). If it is furthermore assumed that $X$ is isotropic in space, then
\begin{equation*}
g(y_1, u_1, y_2, u_2) = g_*\{\|y_1-y_2\|_d, d(u_1, u_2)\}, \qquad y_1,y_2\in\mathbb R^d,\ u_1,u_2\in\mathbb S^k,
\end{equation*}
depends only on $(y_1, y_2)$ through $\|y_1-y_2\|_d$ and on $(u_1, u_2)$ through $d(u_1, u_2)$.

The spatial components of $X$ constitute a usual
\textit{spatial point process} $Y = \{y : (y, u)\in X\}$, which is locally finite, whereas the spherical components constitute a
\textit{point process on the sphere}
$U = \{u: (y, u) \in X\}$ that may be infinite on the compact set $\mathbb S^k$. Let $W\subset\mathbb R^d$ be a bounded Borel set, which we may think of as a window where the spatial components $Y_W=Y\cap W$ are observed. As $X$ is locally finite, the spherical components associated with $Y_W$ constitute a  finite point process $U_W=\{u:(y,u)\in X,y\in W\}$ on $\mathbb S^k$. Let $N=N(W)$ denote the cardinality of $Y_W$.   
To avoid trivial and undesirable cases, we assume that 
$|W|>0$ and that the following inequalities hold:

\begin{equation}\label{e:pos1}
0<\mathrm{E}(N)<\infty
\end{equation} 
and
\begin{equation}\label{e:pos2}
0<\mathrm{E}\{N(N-1)\}<\infty,
\end{equation}
where, by \eqref{eq:defn_intensity}--\eqref{eq:defn_2nd_order_intensity}, 
\begin{align*}
\mathrm{E}(N) = \int_{W\times \mathbb S^k}\rho(y, u)\,\mathrm d\mu(y, u)
\end{align*}
and
\begin{align*}
  \MoveEqLeft \mathrm{E}\{N(N-1)\}
  \\
  &= \int_{W \times \mathbb
    S^{k}}\int_{W\times \mathbb S^{k}}\rho(y_1, u_1)\rho(y_2,
  u_2)g(y_1, u_1, y_2, u_2) \,\mathrm d\mu(y_1, u_1) \,\mathrm
  d\mu(y_2, u_2).
\end{align*}

Note that $Y$ has intensity function $\rho_1$ and pair correlation function $g_1$ given by
\begin{equation}\label{e:rho1}
\rho_1(y)=\int \rho(y,u)\,\mathrm d\nu(u),\qquad y\in \mathbb R^d,
\end{equation}
and
\begin{equation}\label{e:g1}
  \rho_1(y_1)\rho_1(y_2)g_1(y_1,y_2)=
  \iint \rho(y_1, u_1)\rho(y_2, u_2)g(y_1, u_1, y_2, u_2) 
  \,\mathrm d\nu(u_1)\,\mathrm d\nu(u_2)
\end{equation}
for  $y_1,y_2\in \mathbb R^d$,
where we set $g_1(y_1, y_2)=0$ if $\rho_1(y_1)\rho_1(y_2)=0$.
This follows from \eqref{eq:defn_intensity}--\eqref{eq:defn_2nd_order_intensity} and  definitions of the intensity and pair correlation function for spatial point processes \citep[see e.g.][]{MW2004}.
Clearly, if $X$ is stationary in space, then $Y$ is stationary, $\rho_1$ is constant, and $g_1(y_1,y_2)$ is stationary, that is, it depends only on $y_1-y_2$. If in addition $X$ is isotropic in space, then $g_1(y_1,y_2)$ is isotropic, that is, it depends only on $\|y_1-y_2\|_d$. On the other hand if $Y$ is stationary (or isotropic) and the spherical components are independent of $Y$, then $X$ is stationary (or isotropic) in space.

Similarly, using definitions of the intensity and pair correlation function for point processes on the sphere \citep{LBMN2016, MR2016}, $U_W$ has intensity function $\rho_2$ (with respect to $\nu$) and pair correlation function $g_2$ 
(with respect to $\nu\otimes\nu$) given by 
\begin{equation}\label{e:rho2}
\rho_2(u)=\int_W\rho(y,u)\,\mathrm dy,\qquad u\in\mathbb S^k,
\end{equation}
and
\begin{equation}\label{e:g2}
\rho_2(u_1)\rho_2(u_2)g_2(u_1,u_2)=
\int_W \int_W \rho(y_1, u_1)\rho(y_2, u_2)g(y_1, u_1, y_2, u_2) 
\,\mathrm dy_1\,\mathrm dy_2
\end{equation}
for $u_1,u_2\in\mathbb S^k$, where we set $g_2(u_1, u_2)=0$ if $\rho_2(u_1)\rho_2(u_2)=0$. Note that we suppress in the notation that $\rho_2$ and $g_2$ depend on $W$. Obviously, if $X$ is isotropic on the sphere, then $\rho_2$ is constant and $g_2(u_1,u_2)$ is isotropic as it depends only on $d(u_1,u_2)$.

\subsection{Examples}\label{sec:examples}
The following examples introduce the point process models considered in this paper. 

\begin{example}[(Poisson and Cox processes)]\label{ex:PoissonCox}
	First, suppose $X$ is a {\it Poisson process} with a locally integrable intensity function $\rho$. This means that, the count $N(A)=\#X_A$ is Poisson distributed with mean $\int_A\rho(y,u)\,\mathrm d\mu(y, u)$ for any bounded Borel set $A\subset S$ and, conditional on $N(A)$, the points in $X_A$ are independent and identically distributed (IID) with a density proportional to $\rho$ restricted to $A$. Note that $g=1$. Further, $X$ is stationary in space and isotropic on the sphere if and only if $\rho$ is constant, in which case we call $X$ a {\it homogeneous Poisson process} with intensity $\rho$. Furthermore, $Y$ and $U_W$ are Poisson processes, so $g_1=1$ and $g_2=1$.
	
	Second, let $\Lambda=\{\Lambda(y,u):(y,u)\in S\}$ be a non-negative random field so that with probability one $\int_A\Lambda(y,u)\,\mathrm d\mu(y, u)$ is finite for any bounded Borel set $A\subset S$. 
	If $X$ conditioned on $\Lambda$ is a Poisson process with intensity function $\Lambda$, then $X$ is said to be a {\it Cox process driven by $\Lambda$} \citep{Cox1955}. Clearly, the intensity and pair correlation functions of $X$ are 
	\begin{equation}\label{e:cox1}
	\rho(y,u)=\mathrm{E}\{\Lambda(y,u)\}, \qquad y\in\mathbb R^d,\ u\in\mathbb S^k,
	\end{equation}
	and
	\begin{equation*}
	\rho(y_1, u_1)\rho(y_2, u_2)g(y_1, u_1, y_2, u_2)= \mathrm{E}\{\Lambda(y_1,u_1)\Lambda(y_2,u_2)\}, 
	\end{equation*}
	for $y_1,y_2\in\mathbb R^d,\ u_1,u_2\in\mathbb S^k$.
	To separate the intensity function $\rho$ from random effects, it is convenient to work with a so-called residual random field $R = \{ R(y, u) : (y, u) \in S \}$ fulfilling $\Lambda(y,u)=\rho(y,u)R(y,u)$, so $\mathrm{E}\{R(y, u)\} = 1$ \citep[see e.g.][]{moeller:waagepetersen:07, diggle:14}. Then 
	\begin{equation}\label{e:cox2}
	g(y_1, u_1, y_2, u_2)= \mathrm{E}\{R(y_1,u_1)R(y_2,u_2)\},\qquad y_1,y_2\in\mathbb R^d,\ u_1,u_2\in\mathbb S^k,
	\end{equation}
	$\mbox{whenever } \rho(y_1, u_1)\rho(y_2, u_2)>0$. 
\end{example}

Note that projected point processes $Y$ and $U_W$ are Cox processes driven by the random fields $\{ \int_{\mathbb{S}^k} \Lambda(y, u) \, \mathrm d\nu (u) : y \in \mathbb{R}^d\}$ and $\{ \int_{W} \Lambda(y, u) \, \mathrm dy : u \in \mathbb{S}^k\}$, respectively.  Their intensity and pair correlation functions are specified by \eqref{e:rho1}--\eqref{e:g2}.

\begin{example}[(Log Gaussian Cox processes)]\label{ex:LGCP}
	A Cox process $X$ is called a {\it log Gaussian Cox process} 
	\citep[\textit{LGCP};][]{moeller:syversveen:waagepetersen:98}
	if the residual random field is of the form $R = \exp(Z)$, where $Z$ is a Gaussian random field (GRF) with mean function $\mu(y, u) = - c(y, u, y, u)/2$, where $c$ is the covariance function of $Z$. Note that $X$ has pair correlation function 
	\begin{equation}\label{e:lgcp2}
	g(y_1, u_1, y_2, u_2) = \exp \left[ c\{ (y_1, u_1), (y_2, u_2) \} \right], \qquad y_1, y_2 \in \mathbb R^d, \, u_1, u_2 \in \mathbb S^k.
	\end{equation}
\end{example}

\begin{example}[(Marked point processes)]\label{ex:IndMark}	
	It is sometimes useful to view $X$ as a marked point process \citep[see e.g.][]{daley:vere-jones:03,illian:etal:08}, where the spatial components are treated as the ground process and the spherical components as marks. Often it is of interest to test the hypothesis $H_0$ that the marks are IID and independent of the ground process $Y$. 
	Under $H_0$, with each mark following a density $p$ with respect to $\nu$, the intensity is
	\[
	\rho(y, u) = \rho_1(y)p(u), \qquad y \in \mathbb R^d, \, u \in \mathbb S^k,
	\]
	and the pair correlation function 
	\[	
	g(y_1, u_1, y_2, u_2) = g_1(y_1, y_2), \qquad y_1,y_2 \in \mathbb R^d, \, u_1, u_2 \in \mathbb S^k,	
	\]
	does not depend on $(u_1, u_2)$.
		
	In some situations, it may be more natural to look at it conversely, that is, treating  $U_W$ as the ground process and $Y_W$ as marks. Then similar results for $\rho$ and $g$ may be established by interchanging the roles of points and marks. 
\end{example}

\begin{example}[(Independently marked determinantal point processes)]\label{ex:DPP}
  Considering a space-sphere point process $X$ as a marked point
  process that fulfils the hypothesis $H_0$ given in
  Example~\ref{ex:IndMark}, we may let the ground process $Y$ be
  distributed according to any point process model of our choice
  regardless of the marks $U$. For instance, in case of repulsion
  between the points in $Y$, a \textit{determinantal point process
    (DPP)} may be of interest because of its attractive properties
  \citep[see][and the references therein]{LMR2015}.  Briefly, a DPP is
  defined by a so-called kernel
  $C: \mathbb R^d \times \mathbb R^d \to \mathbb C$, which we assume
  is a complex covariance function, that is, $C$ is positive
  semi-definite and Hermitian. Furthermore, let $\rho_1^{(n)}$ denote
  the $n$th order joint intensity function of $Y$, that is,
  $\rho_1\psup{1} = \rho_1$ is the intensity and
  $\rho_1\psup{2}(y_1, y_2) = \rho_1(y_1)\rho_1(y_2)g_1(y_1, y_2)$ for
  $y_1, y_2 \in \mathbb R^d$, while we refer to \cite{LMR2015} for the
  general definition of $\rho_1\psup{n}$ which is an extension of
  \eqref{e:rho1}--\eqref{e:g1}.  If for all $n =1, 2, \dots$,
	\[
          \rho_1\psup{n}(y_1, \dots, y_n) = \det\{C(y_i, y_j)\}_{i,j =
            1, \dots, n}, \qquad y_1, \dots, y_n \in \mathbb{R}^d,
	\]
	where $\det\{C(y_i, y_j)\}_{i,j = 1, \dots, n}$ is the determinant of the $n\times n$ matrix with $(i, j)$-entry $C(y_i, y_j)$, we call $Y$ a DPP with kernel $C$ and refer to $X$ as an \textit{independently marked DPP}.  
	It follows that $Y$ has intensity function $\rho(y)=C(y, y)$ and pair correlation function
	\begin{align*}
	g_1(y_1, y_2) = 1 -  | R(y_1, y_2) |^2, \qquad y_1, y_2 \in \mathbb{R}^d, 
	\end{align*}
	whenever $\rho(y_1)\rho(y_2)>0$, where $R(y_1, y_2) = C(y_1, y_2)/\sqrt{C(y_1, y_1)C(y_2, y_2)}$ is the correlation function corresponding to $C$ and $|z|$ denotes the modulus of $z\in \mathbb C$. 
	
	Alternatively, we may look at a DPP on the sphere \citep{MNPR2018}, that is, modelling $U_W$ as a DPP while considering $Y_W$ as the marks and impose the conditions of IID marks independent of $U_W$.
\end{example}

\section{The space-sphere $K$-function}\label{s:SOIRS}

\subsection{Definition}\label{sec:SOIRSdef}

When \eqref{e:gSOIRS} holds we say that the space-sphere point process $X$ is \textit{second order intensity-reweighted stationary} (SOIRS) and
define the \textit{space-sphere $K$-function} by 
\begin{align}\label{e:Krs}
K(r, s) = 
\int_{\|y\|_d\leq r, \, d(u, e) \leq s }g_0\{y, d(u, e)\}\,\mathrm d\mu(y, u),\qquad  r \geq 0, \, 0\leq s \leq \pi, 
\end{align}
where $e \in \mathbb{S}^k$ is an arbitrary reference direction. This
definition does not depend on the choice of $e$, as the integrand only
depends on $u\in \mathbb{S}^k$ through its geodesic distance to $e$
and $\nu(\cdot)$ is a rotation invariant measure. For example, we may
let $e=(0,\ldots,0,1)\in\mathbb S^k$ be the ``North Pole''.

Let $\sigma_k = \nu(\mathbb{S}^k)=2\pi^{(k+1)/2}/\Gamma\{(k+1)/2\}$ denote the surface measure of $\mathbb S^k$. For any Borel set $B\subset\mathbb R^d$ with $0<|B|<\infty$, we easily obtain from \eqref{eq:defn_2nd_order_intensity} and \eqref{e:Krs} that
\begin{align}
  \MoveEqLeft[1] K(r, s) \nonumber
  \\
  &= \frac{1}{|B|\sigma_k}\iint_{y_1 \in B, \, \|y_1 - y_2\|_d \leq r,
    \, d(u_1, u_2) \leq s} g_0\{y_1-y_2,
  d(u_1, u_2)\}\,\mathrm d\mu(y_1, u_1)\,\mathrm d\mu(y_2,
  u_2)\nonumber
  \\[0.3em]
  &= \frac{1}{|B|\sigma_k}\mathrm{E} \Bigl[\sum^{\neq}_{(y_i, u_i),
    (y_j, u_j) \in X}\frac{\mathbb{I}\{y_i \in B, \, \| y_i-y_j\|_d
    \leq r, \, d(u_i, u_j) \leq s\}}{\rho(y_i, u_i)\rho(y_j,
    u_j)}\Bigr] \label{e:B}
\end{align}
for $r\geq 0, \, 0 \leq s \leq \pi$, where $\mathbb{I}(\cdot)$ denotes the indicator function.
The relation given by \eqref{e:B} along with the requirement that the expression in \eqref{e:B} does not depend on the choice of $B$ could alternatively have been used as a more general definition of the space-sphere $K$-function. Such a definition is in agreement with the one used in  \cite{BMW2000} for SOIRS of a spatial point process. It is straightforward to show that \eqref{e:B} does not depend on $B$ when $X$ is stationary in space.

For $r, s >0$ and $(y_1, u_1), (y_2, u_2) \in S$, we say that $(y_1, u_1)$ and $(y_2, u_2)$ are $(r, s)$-close neighbours if $\| y_1 - y_2\|_d \leq r$ and $d(u_1, u_2) \leq s$. If $X$ is stationary in space and isotropic on the sphere, then \eqref{e:B} shows that 
$\rho K(r, s)$ can be interpreted as the expected number of further $(r, s)$-close neighbours in $X$ of a typical point in $X$. More formally, this interpretation relates to the reduced Palm distribution \citep{daley:vere-jones:03}.

Some literature treating marked point processes discuss the so-called
\textit{mark-weighted $K$-function} \citep[see e.g.][]{illian:etal:08,
  Koubek2016}, which to some extent resembles the space-sphere
$K$-function in a marked point process context; both are cumulative
second order summary functions that consider points as well as
marks. However, the mark-weighted $K$-function has an emphasis on the
marked point process setup (and considers e.g.\ $\rho_1$ rather than
$\rho$), whereas the space-sphere $K$-function is constructed in such
a way that it is an analogue to the planar/spherical $K$-function for
space-sphere point processes.

\begin{myexample}[Example~\ref{ex:PoissonCox} continued (Poisson and Cox processes)]
A Poisson process is clearly SOIRS and $K(r, s)$ is simply the product of the volume of a $d$-dimensional ball with radius $r$ and the surface area of a spherical cap given by $\{u\in\mathbb S^k: d(u, e) \leq s\}$ for an arbitrary $e \in \mathbb S^k$ \citep[see][for formulas of this area]{Li2011}. Thus, for $r\geq 0$, the space-sphere $K$-function is
	\begin{equation*}
	K_{Pois}(r, s) = 
	\begin{cases}
	\frac{r^d\pi^{(d+k+1)/2}}{\Gamma(1 + d/2)\Gamma\{(k+1)/2\}} I_{\sin^2(s)}\left(\frac{k}{2}, \frac{1}{2}\right), \qquad 0\leq s \leq \frac{\pi}{2},\\[0.5em]
	\frac{r^d\pi^{(d+k+1)/2}}{\Gamma(1 + d/2)\Gamma\{(k
          +1)/2\}}\{2 - I_{\sin^2(\pi - s)}(\frac{k}{2},
            \frac{1}{2})\}, \qquad \frac{\pi}{2} < s \leq
        \pi,
	\end{cases}
	\end{equation*}
	where $I_x(a, b)$ is the regularized incomplete beta function. 
	In particular, if $k = 2$, 
	\[
	\begin{rcases}
	I_{\sin^2(s)}(\frac{k}{2}, \frac{1}{2}), & 0 \leq s \leq \frac{\pi}{2} \\
	2 - I_{\sin^2(\pi - s)}(\frac{k}{2}, \frac{1}{2}),
        & \frac{\pi}{2}< s \leq \pi
	\end{rcases}
	= 
	1 - \cos(s).
	\]
	
	If the residual random field $R$ in \eqref{e:cox2} is invariant under translations in $\mathbb R^d$ and under rotations on $\mathbb S^k$, then the associated Cox process is SOIRS. The evaluation of $g$ (and thus $K$) depends on the particular model of $R$ as exemplified in Example~\ref{ex:LGCP} below and in Section~\ref{s:discussion}.
\end{myexample}

\begin{myexample}[Example~\ref{ex:LGCP} continued (LGCPs)]
	Suppose that the distribution of $R$ is invariant under translations in $\mathbb R^d$ and under rotations on $\mathbb S^k$, and recall that $R$ is required to have unit mean. Then the underlying GRF $Z$ has a covariance function of the form 
	\[
	c(y_1,u_1,y_2,u_2)=c_0\{y_1-y_2,d(u_1,u_2)\},\qquad y_1,y_2\in\mathbb R^d,\ u_1,u_2\in\mathbb S^k,
	\]
	and $\mathrm{E}Z(y,u)=-\sigma^2/2$ for all $y \in \mathbb R^d$ and $u\in\mathbb S^k$, where $\sigma^2=c_0(0,0)$ is the variance. It then follows from \eqref{e:lgcp2} that $X$ is SOIRS with
	\begin{equation}\label{e:lgcp3}
	g_0(y,s)=\exp\{c_0(y,s)\},\qquad y\in\mathbb R^d,\ 0\le s\le\pi.
	\end{equation}
\end{myexample}

\section{Separability}\label{s:separability}
\subsection{First order separability}\label{s:firstordersep}
We call the space-sphere point process $X$ \textit{first order separable} if there exist non-negative Borel functions $f_1$ and $f_2$ such that
\begin{align*}
\rho(y, u) = f_1(y)f_2(u),\qquad y\in \mathbb R^d,\ u\in\mathbb S^k.
\end{align*}
By \eqref{e:pos1}, \eqref{e:rho1}, and \eqref{e:rho2} this is equivalent to 
\begin{equation}\label{e:sep1}
\rho(y, u) = \rho_1(y)\rho_2(u)/{\mathrm{E}(N)},\qquad y\in \mathbb R^d,\ u\in\mathbb S^k,
\end{equation}
recalling that $\rho_2$ and $N$ depend on $W$, but $\rho_2/\mathrm{E}N$ does not depend on the choice of $W$.
Then, in a marked point process setup where the spherical components are treated as marks,  $\rho_2(\cdot)/{\mathrm{E}(N)}$ is the density of the mark distribution.
First order separability was seen in Example~\ref{ex:IndMark} to be fulfilled under the assumption of IID marks independent of the ground process. Moreover, any homogeneous space-sphere point process is clearly first order separable. 
In practice, first order separability is a working hypothesis which may be hard to check. 

\subsection{Second order separability}\label{s:secondordersep}
If there exist Borel functions $k_1$ and $k_2$ such that
\begin{align*}
g(y_1, u_1, y_2, u_2) = k_1(y_1, y_2)k_2(u_1, u_2),\qquad y_1,y_2\in \mathbb R^d,\ u_1,u_2\in\mathbb S^k,
\end{align*}
we call $X$ \textit{second order separable}. Assuming first order separability, it follows by \eqref{e:pos2}, \eqref{e:g1}, \eqref{e:g2}, and \eqref{e:sep1} that second order separability is equivalent to
\begin{equation}\label{e:sep2}
g(y_1, u_1, y_2, u_2) = \beta g_1(y_1, y_2)g_2(u_1, u_2),\qquad y_1,y_2\in \mathbb R^d,\ u_1,u_2\in\mathbb S^k,
\end{equation}
where 
\[
\beta={\mathrm{E}(N)^2}/{{\mathrm{E}\{N(N-1)\}}}
\]
and noting that $\beta$ and $g_2$ depend on $W$, but $\beta g_2$ does not depend on the choice of $W$.
The value of $\beta$ may be of interest: for a Poisson Process, $\beta = 1$; for a Cox process, $\mathrm{var}(N) \geq \mathrm{E}(N)$ \citep[see e.g.][]{MW2004}, so $\beta \leq 1$; for an independently marked DPP, $\beta \geq 1$ \citep{LMR2015}.

\begin{myexample}[Example~\ref{ex:PoissonCox} continued (Poisson and Cox processes)]
  Clearly, when $X$ is a Poisson process, it is second order
  separable.  Assume instead that $X$ is a Cox process and the
  residual random field is separable, that is, $R(y,u)=R_1(y)R_2(u)$,
  where $R_1=\{R_1(y):y\in\mathbb R^d\}$ and
  $R_2=\{R_2(u):u\in\mathbb S^k\}$ are independent random fields.
  Then, by \eqref{e:cox2}, $X$ is second order separable and
	\begin{align*}
          \MoveEqLeft g(y_1, u_1, y_2, u_2)
          \\
          &=
          \mathrm{E}\{R_1(y_1)R_1(y_2)\}\mathrm{E}\{R_2(u_1)R_2(u_2)\},
          \qquad y_1, y_2 \in \mathbb R^d,\ u_1, u_2 \in \mathbb
          S^k.
	\end{align*} 
\end{myexample}

\begin{myexample}[Example~\ref{ex:LGCP} continued (LGCPs)]
  If $X$ is a LGCP driven by $\Lambda(y,u)=\rho(y,u)\cdot\exp\{Z(y,u)\}$,
  second order separability is implied if $Z_1=\log R_1$ and
  $Z_2=\log R_2$ are independent GRFs so that
  $Z(y,u)=Z_1(y)+Z_2(u)$. Then, by the imposed invariance properties
  of the distribution of the residual random field, $Z_1$ must be
  stationary with a stationary covariance function
  $c_1(y_1,y_2)=c_{01}(y_1-y_2)$ and mean $-c_{01}(0)/2$, and $Z_2$
  must be isotropic with an isotropic covariance function
  $c_2(u_1,u_2)=c_{02}\{d(u_1,u_2)\}$ and mean $-c_{02}(0)/2$.
  Consequently, in \eqref{e:lgcp3}, $c_0(y,s)=c_{01}(y)+c_{02}(s)$ for
  $y \in \mathbb R^d$ and $0 \leq s \leq \pi$.
\end{myexample}

\begin{myexample}[Example~\ref{ex:IndMark} continued (marked point processes)]
  Consider the space-sphere point process $X$ as a marked point
  process with marks in $\mathbb S^k$.  As previously seen, first and
  second order separability is fulfilled under the assumption of IID
  marks independent of the ground process, but we may in fact work
  with weaker conditions to ensure the separability properties as
  follows. Assume that each mark is independent of the ground process
  $Y$ and the marks are identically distributed following a density
  function $p$ with respect to $\nu$.  Then the first order
  separability condition \eqref{e:sep1} is satisfied with
  $\rho_2(u) = \mathrm{E}(N)p(u)$ for $u \in \mathbb S^k$.  In addition,
  assuming the conditional distribution of the marks given $Y$ is such
  that any pair of marks is independent of $Y$ and follows the same
  joint density $q(\cdot,\cdot)$ with respect to $\nu\otimes\nu$, it
  is easily seen that the second order separability condition
  \eqref{e:sep2} is satisfied with
  \[
    g_2(u_1,u_2)=\frac{q(u_1,u_2)}{\beta p(u_1)p(u_2)}, \qquad
    u_1,u_2\in\mathbb S^k,
  \]
  whenever $\rho_2(u_1)\rho_2(u_2)>0$.  If we also have pairwise
  independence between the marks, that is, $q(u_1,u_2)=p(u_1)p(u_2)$,
  then the pair correlation function
  $g(y_1,u_1,\allowbreak y_2,u_2)=g_1(y_1,y_2)$ does not depend on
  $(u_1, u_2)$ and $g_2(u_1, u_2) = 1/\beta$ is constant. Note that
  this implies $g_2 \leq 1$ for an independently marked DPP, reflecting that even when the marks are drawn independently of $Y$ the behaviour of the points implicitly affects the marks as the number of points is equal to the number of marks.

  Again, the roles of points and marks may be switched resulting in
  statements analogue to those above.
\end{myexample}

\subsection{Assuming both SOIRS and first and second order separability}\label{s:SOIRS_sep}
Suppose that $X$ is both SOIRS and first and second order separable. Then the space-sphere $K$-function can be factorized as follows. Note that $Y$ and $U_W$ are SOIRS since there by \eqref{e:gSOIRS}, \eqref{e:g1}, \eqref{e:g2}, and \eqref{e:sep1} 
exist Borel functions $g_{01}$ and $g_{02}$ such that 
\begin{equation}\label{e:g01}
  \begin{aligned}
    g_1(y_1, y_2) &= g_{01}(y_1 - y_2)
    \\
    &=\iint\frac{\rho_2(u_1)}{\mathrm{E}(N)}\frac{\rho_2(u_2)}{\mathrm{E}(N)}
    g_0\{y_1-y_2,d(u_1,u_2)\}\,\mathrm d\nu(u_1)\,\mathrm d\nu(u_2)
  \end{aligned}
\end{equation}
for $y_1,y_2\in \mathbb R^d\mbox{ with }\rho_1(y_1)\rho_1(y_2)>0$, 
and  
\begin{equation}\label{e:g02}
  \begin{aligned}
    g_2(u_1, u_2)& = g_{02}\{d(u_1, u_2)\}
    \\
    &=\int_{W}\int_{W}\frac{\rho_1(y_1)}{\mathrm{E}(N)}\frac{\rho_1(y_2)}{\mathrm{E}(N)}
    g_0\{y_1-y_2,d(u_1,u_2)\}\,\mathrm dy_1\,\mathrm
    dy_2
  \end{aligned}
\end{equation}
for $u_1,u_2\in\mathbb S^k\mbox{ with }\rho_2(u_1)\rho_2(u_2)>0$. 
Hence, the inhomogeneous $K$-function for the spatial components in $Y$ \citep[introduced in][]{BMW2000} is
\begin{align*}
K_1(r) = \int_{\|y\|_d\leq r}g_{01}(y)\,\mathrm dy,\qquad r \geq 0,
\end{align*}
and the inhomogeneous $K$-function for the spherical components in $U_W$ \citep[introduced in][]{LBMN2016, MR2016} is
\begin{align*}
K_2(s) = \int_{d(u, e) \leq s} g_{02}\{d(u, e)\}\,\mathrm d\nu(u),\qquad 0 \leq s \leq \pi, 
\end{align*}
where $e \in \mathbb{S}^k$ is arbitrary. 
Combining \eqref{e:Krs} and \eqref{e:sep2}--\eqref{e:g02},
we obtain   
\begin{align*}
K(r, s) = \beta K_1(r)K_2(s),\qquad r \geq 0, \, 0\leq s \leq \pi.
\end{align*}
Note that, if $X$ is a first order separable Poisson process, then $D(r, s) = K(r, s) - K_1(r)K_2(s)$ is 0, and an estimate of $D$ may also be used as a functional summary statistic when testing a Poisson hypothesis.

\section{Estimation of $K$-functions}\label{s:est}
In this section,  we assume for specificity that the observation window is $W \times \mathbb{S}^k$, where $W\subset \mathbb{R}^d$ is a bounded Borel set, and a realisation $X_{W\times\mathbb S^k}=x_{W\times\mathbb S^k}$ is observed; in Section~\ref{s:discussion}, we discuss other cases of observation windows. We let $Y_W=y_W$ and $U_W=u_W$ be the corresponding sets of observed spatial and spherical components. 

First, assume that $\rho_1$ and $\rho_2$ are known. Following \cite{BMW2000}, we estimate $K_1$ by 
\begin{align}\label{e:estK1}
\hat{K}_1(r)  = \sum^{\neq}_{y_i, y_j \in y_W}\frac{\mathbb{I}(\| y_i - y_j\|_d\leq r)}{w_1(y_i, y_j)\rho_1(y_i)\rho_1(y_j)}, \qquad r \geq 0,
\end{align}
where $w_1$ is an edge correction factor on $\mathbb{R}^d$. If we let
$w_1(y_i, y_j) = |W\cap W_{y_i - y_j}|$ be the translation correction
factor \citep{Ohser1983}, where $W_y = \{y + z : z\in W\}$ denotes the
translation of $W$ by $y \in \mathbb{R}^d$, then $\hat{K}_1$ is an
unbiased estimate of $K_1$ \citep[see~e.g.~Lemma~4.2 in][]{MW2004}.
For $d = 1$, we may instead use the temporal edge correction factor
with $w_1(y_i, y_j) = |W|$ if $[y_i - y_j, y_i + y_j] \subseteq W$ and
$w_1(y_i, y_j) = |W|/2$ otherwise \citep{Diggle1995, Moeller2012}.
Moreover, for estimation of $K_2$, we use the unbiased estimate
\begin{align}\label{e:estK2}
\hat{K}_2(s)  = \frac{1}{\sigma_k}\sum^{\neq}_{u_i, u_j \in u_W}\frac{\mathbb{I}\{d(u_i, u_j) \leq s\}}{\rho_2(u_i)\rho_2(u_j)}, \qquad 0 \leq s \leq \pi, 
\end{align}
cf.\ \cite{LBMN2016} and \cite{MR2016}. A natural extension of the above estimates gives the following estimate of $K$:
\begin{align}\label{e:estK}
\hat{K}(r, s)  = \frac{1}{\sigma_k}\sum^{\neq}_{(y_i, u_i), (y_j, u_j) \in x_{W\times\mathbb S^k}}\frac{\mathbb{I}\{\| y_i - y_j\|_d\leq r, \, d(u_i, u_j) \leq s\}}{w_1(y_i, y_j)\rho(y_i, u_i)\rho(y_j, u_j)}
\end{align}
for $r \geq 0, \, 0 \leq s \leq \pi$.
This is straightforwardly seen to be an unbiased estimate when $w_1$ is the translation correction factor.

Second, in practice we need to replace $\rho_1$ in \eqref{e:estK1}, $\rho_2$ in \eqref{e:estK2}, and $\rho$ in \eqref{e:estK} by estimates, as exemplified in Section~\ref{s:appl}. This may introduce a bias.

\section{Data examples}\label{s:appl}

\subsection{Fireball locations over time}\label{sec:fireball}
\begin{figure}
	\centering
	\includegraphics[width=.4\textwidth]{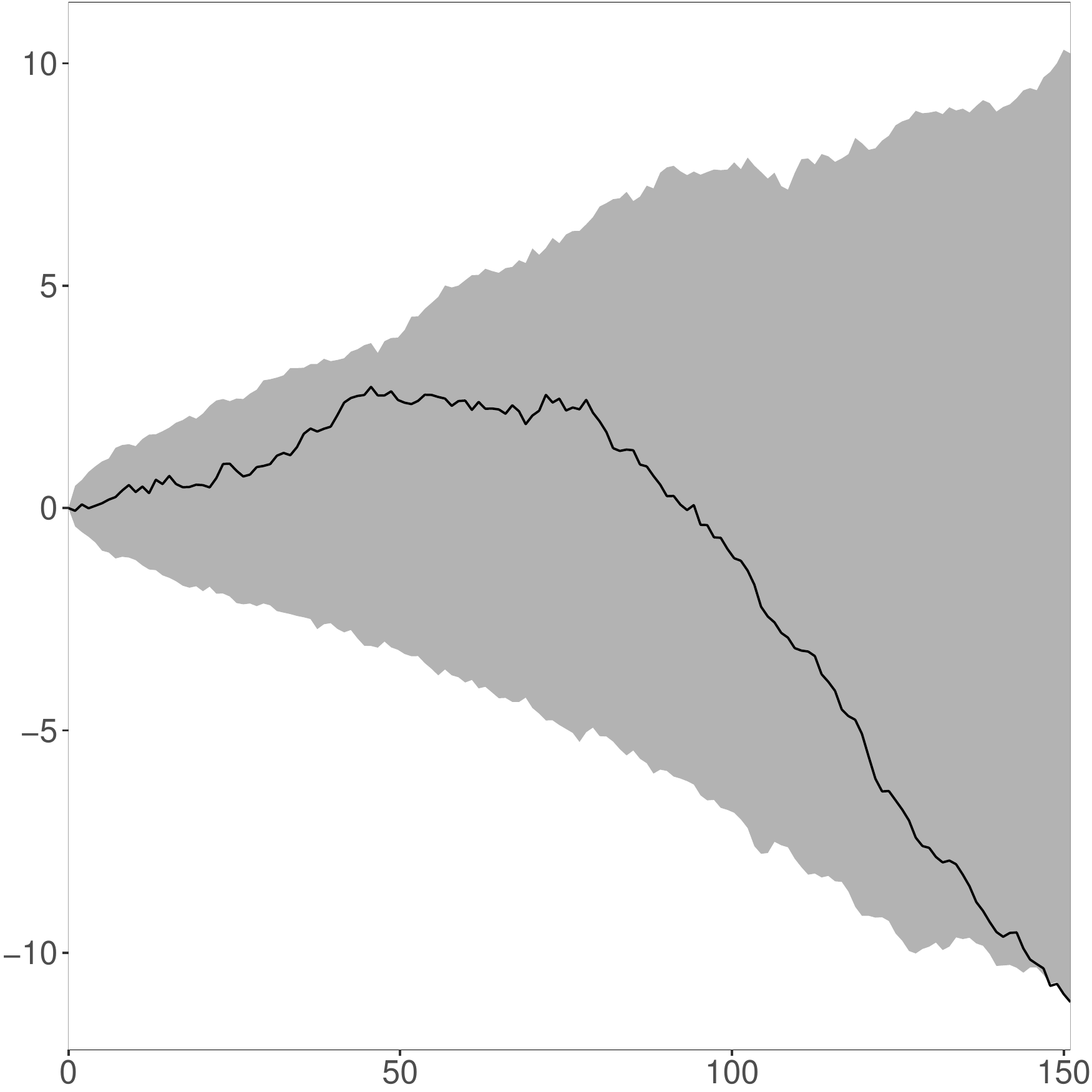}\hspace{5mm}
	\includegraphics[width=.4\textwidth]{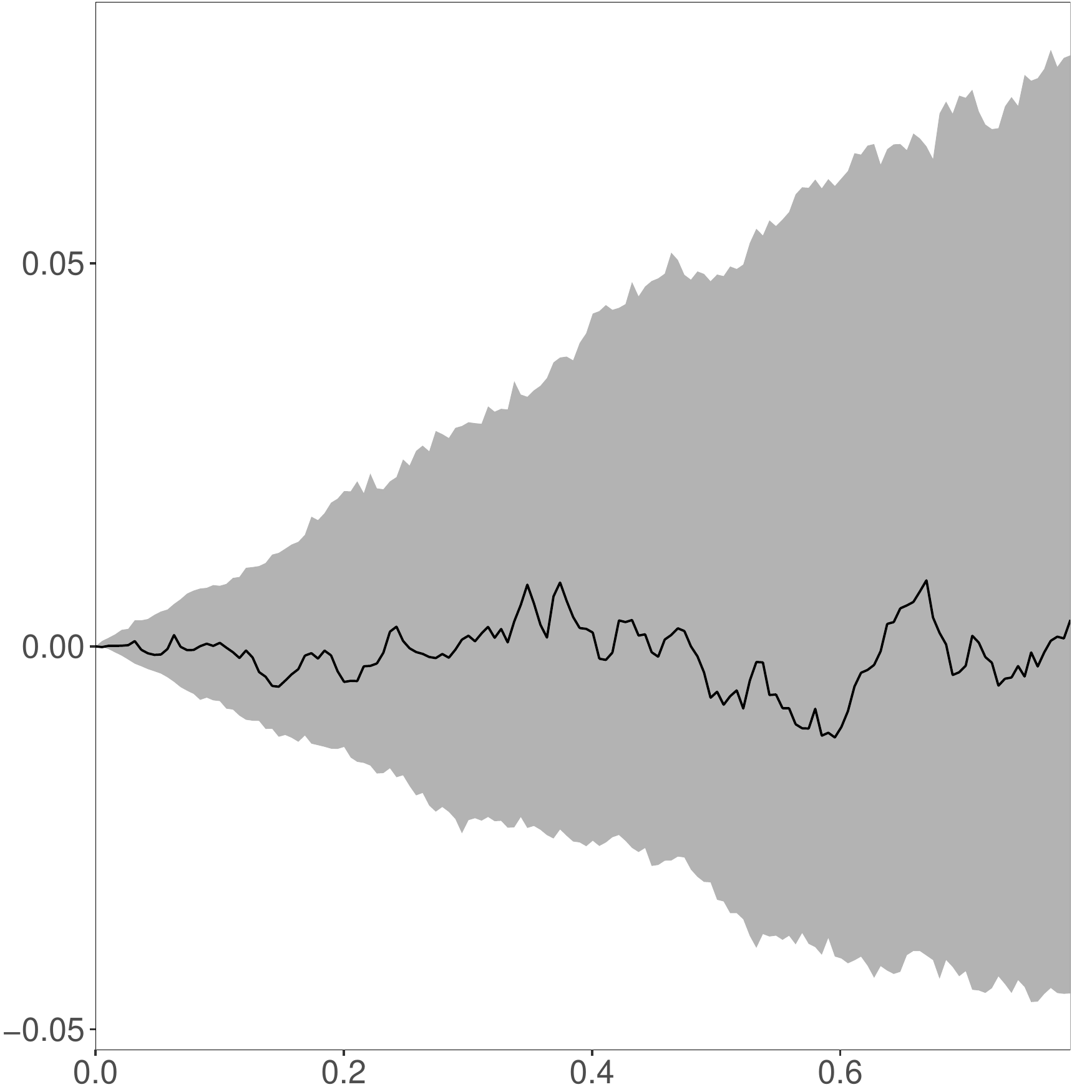}
	\caption{Left: $\hat{K}_{1}(r) - 2r$ for the fireball event times (solid curve) along with a 95\% global rank envelope (grey area) under a homogeneous Poisson  model on the time interval for the observed events. Right: $\hat{K}_{2}(s) - 2\pi \{1 - \cos(s)\}$ for the fireball locations (solid curve) along with a 95\% global rank envelope (grey area) under a homogeneous Poisson model on $\mathbb{S}^{2}$.}
	\label{fig:MarginalKfun}
\end{figure}

\begin{figure}
	\centering
	\includegraphics[width=.49\textwidth]{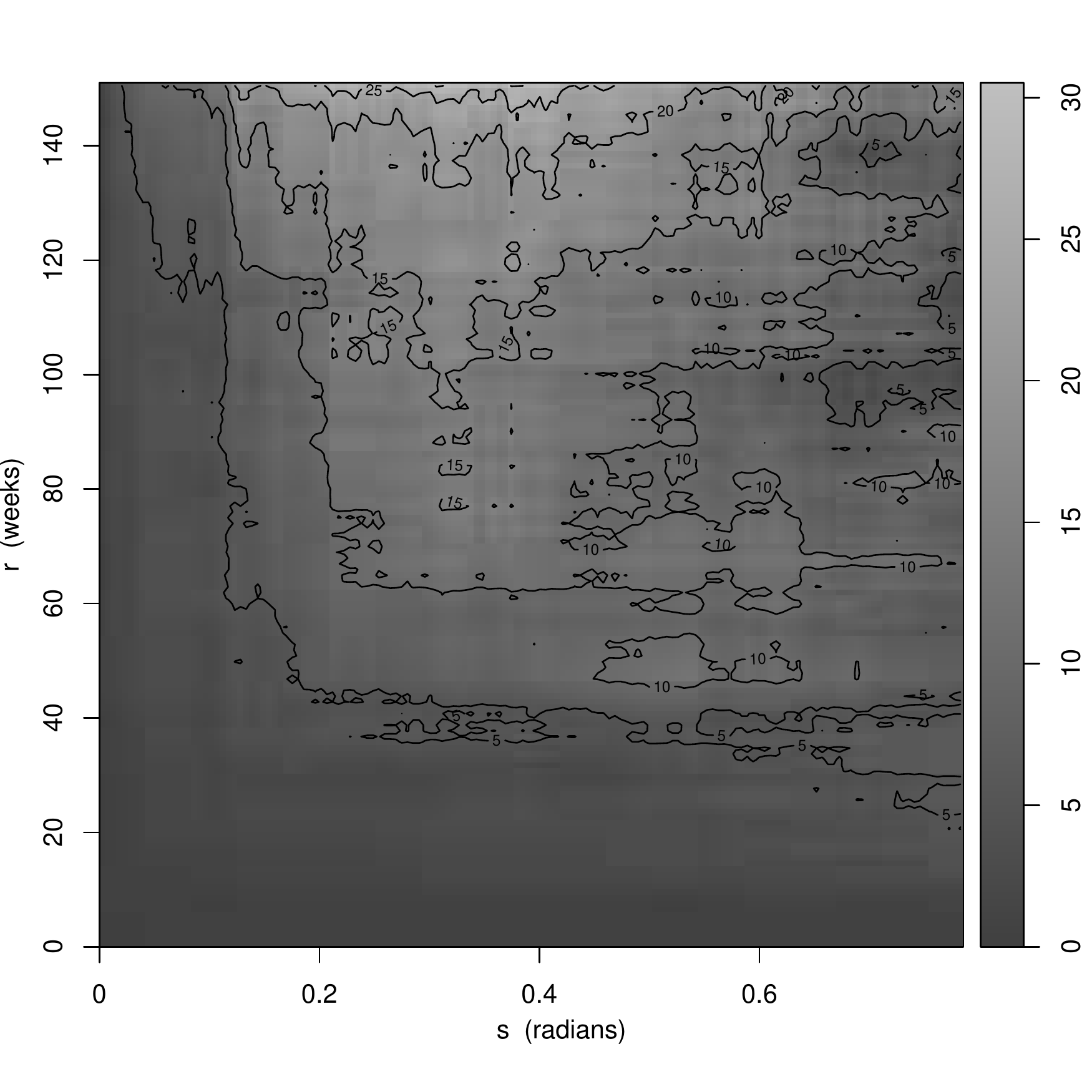}\hspace{1mm}
	\includegraphics[width=.49\textwidth]{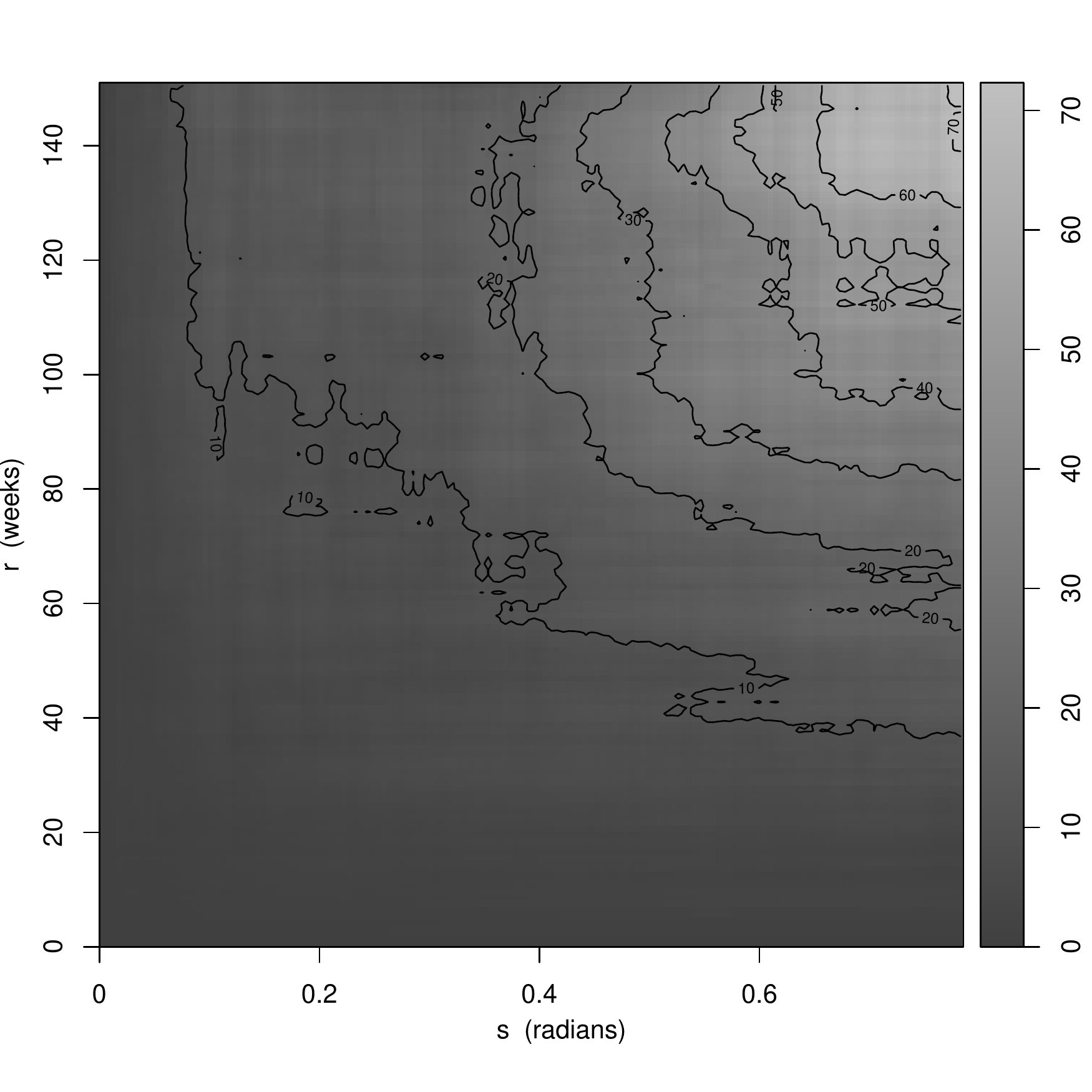}
	\caption{Difference between $\hat{K}(r, s)$ for the observed fireball locations over time and the lower ($\hat{K}_\mathrm{low}$) or upper ($\hat{K}_\mathrm{upp}$) 95\% global rank envelope under a homogeneous Poisson model on $W \times \mathbb{S}^2$. Left: $\hat{K}(r, s) - \hat{K}_\mathrm{low}(r, s)$.  Right: $\hat{K}_\mathrm{upp}(r, s)-\hat{K}(r, s)$.}
	\label{fig:KEnvelope}
\end{figure}

Figure~\ref{fig:MapFireballsSphere} shows the time and location of $n = 344$ fireballs observed over a time period from \texttt{2005-01-01 03:44:09} to \texttt{2016-08-12 23:59:59} corresponding to a time frame $W$ of about 606 weeks. The data can be recovered at \url{http://neo.jpl.nasa.gov/fireballs/}  using these time stamps. 
Figure~\ref{fig:MapFireballsSphere} reveals no inhomogeneity of neither fireball locations or event times. Therefore we assumed first order homogeneity, and used the following unbiased estimates for the intensities: 
\begin{align*}
\hat{\rho}_1 = n/|W| = 0.57, \  \hat{\rho}_2 = n/(4\pi)= 27.37, \ \hat{\rho} = n/(4\pi |W|)= 0.05.
\end{align*}
Then $\hat{K}_1$, $\hat{K}_2$, and $\hat{K}$ (with $w_1$ in \eqref{e:estK1} and \eqref{e:estK} equal to the temporal edge correction factor) were used as test functions in three different global rank envelope tests for testing whether fireball event times, locations, and locations over time each could be described by a homogeneous Poisson model with estimated intensity $\hat{\rho}_1$, $\hat{\rho}_2$, and $\hat{\rho}$, respectively. Appendix A provides a brief account on global rank envelope tests; see also \cite{Myllymaki2017}.
Under each of the three fitted Poisson processes and using 2499 simulations  \citep[as recommended in][]{Myllymaki2017}, we obtained $p$-intervals of $(0.028, 0.040)$ for the event times, $(0.908,  0.908)$ for the locations, and  $(0.445, 0.516)$ for the locations over time. The associated $95 \%$ global rank envelopes for $\hat{K}_1$ and $\hat{K}_2$ are shown in Figure~\ref{fig:MarginalKfun}, and the difference between $\hat{K}$ and the upper and lower $95\%$ global rank envelope is shown in Figure~\ref{fig:KEnvelope}. Since $\hat{K}_2$ and $\hat{K}$ stay inside the $95\%$ global rank envelopes for the considered distances on $\mathbb{S}^k$ and $\mathbb{R}\times \mathbb{S}^k$, there is no evidence against a homogeneous Poisson model for neither locations or locations over time. On the other hand, with a conservative $p$-value of $4\%$, the global rank envelope test based on $\hat{K}_1$ indicates that a homogeneous Poisson model for the event times is not appropriate. 
However, the observed test function $\hat{K}_1(r)$ falls only outside the envelope in Figure~\ref{fig:MarginalKfun} for large values of $r$. Thus, choosing a slightly smaller interval of $r$-values would lead to a different conclusion.

As an alternative to the space-sphere $K$-function, we considered the summary function $D(r,s)$ which in case of a Poisson process is $0$. Estimating $D$ by $\hat{D}(r,s) = \hat{K}(r,s) - \hat{K}_{1}(r)\hat{K}_{2}(s)$, we performed a global rank envelope test with $D$ as test function. The resulting test gave a $p$-interval of $(0.537, 0.564)$ which is similar to the one obtained using $\hat{K}$ as test function.

\subsection{Location and orientation of pyramidal neurons}\label{s:neuron}
We now return to the space-sphere point pattern concerning location
and orientation of pyramidal neurons described in Section~\ref{intro},
which is a data set collected by Ali H.\ Rafati, a biomedical and
clinical scientist. The point pattern is observed on
$W \times \mathbb{S}^2$, where $W \subset \mathbb{R}^3$ is the
rectangular box shown in Figure~\ref{fig:brain_cells} with side
lengths \SI{492.7}{\micro\meter}, \SI{132.0}{\micro\meter}, and
\SI{407.7}{\micro\meter}. Due to the way data was collected, 46
neurons in the original dataset had an orientation/unit vector lying
exactly in the $x$-$z$ plane, meaning that the orientations of these
46 neurons are located only on a great circle of~$\mathbb{S}^2$. To
keep the analysis simple we disregarded these neurons, resulting in a
dataset consisting of $n = 504$ neurons. Below we initially discuss an
appropriate parametric forms of the intensity for the locations and
orientations and how to estimate the intensity parameters. Then we
investigate whether the orientations and locations can be described by
a Poisson model with the proposed intensity, where we first consider
the data as two separate point patterns (a spatial point pattern
describing the locations and a spherical point pattern describing the
orientations) and next as a space-sphere point pattern.

\begin{figure}[!ht]
	\centering
	\includegraphics[width=.55\textwidth, angle=90]{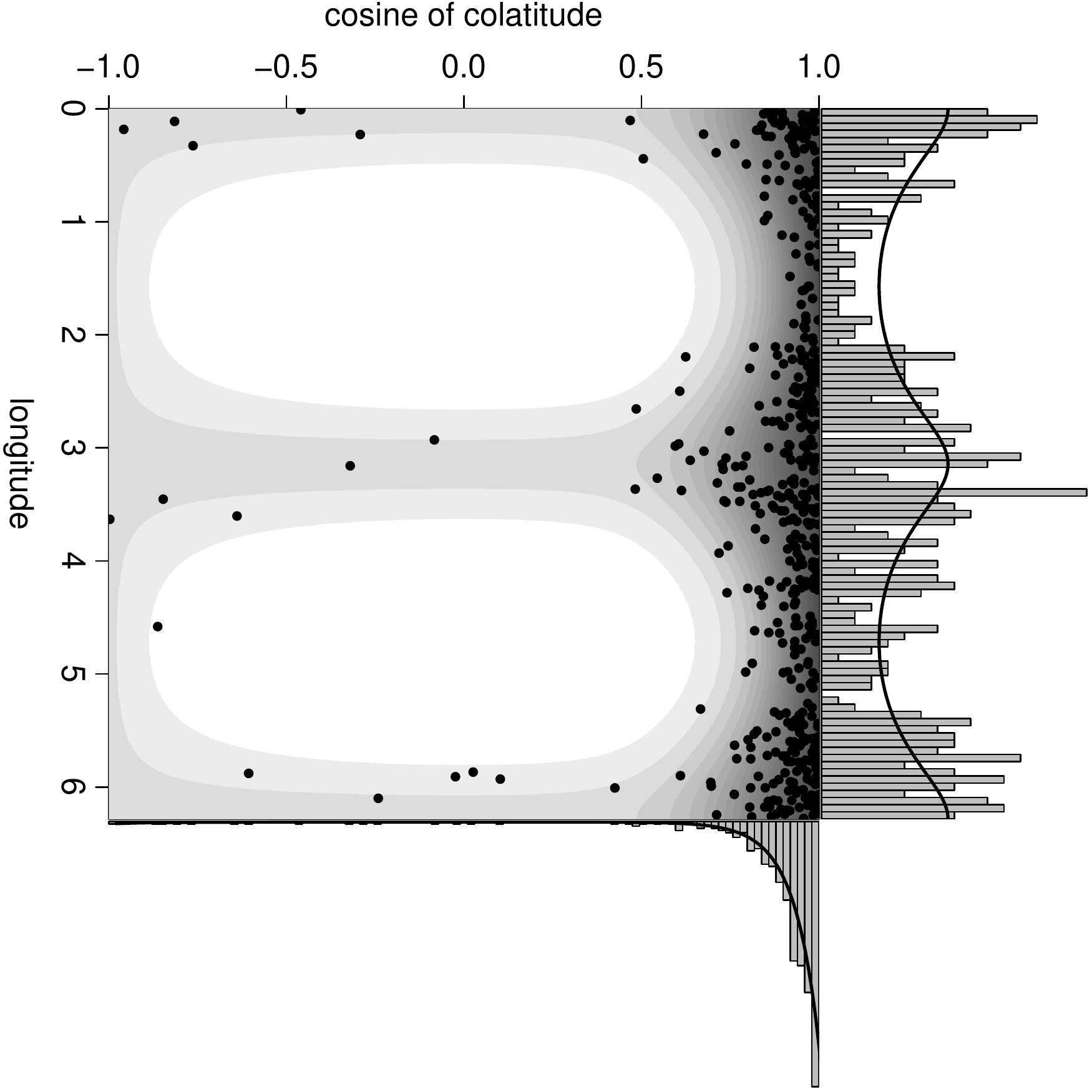}
	\includegraphics[width=.25\textwidth]{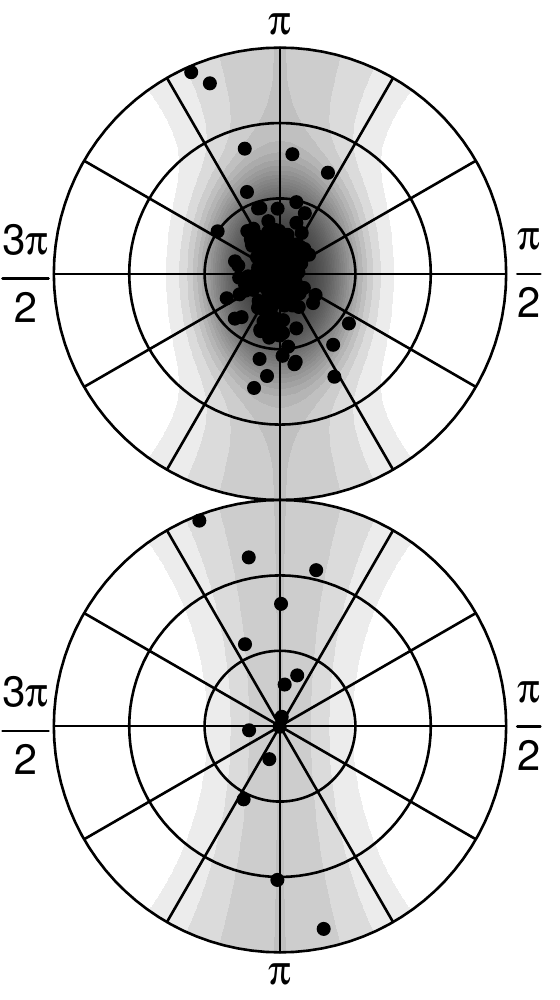}
	\includegraphics[width=.65\textwidth, trim = {0cm 7cm 0cm 7cm}, clip]{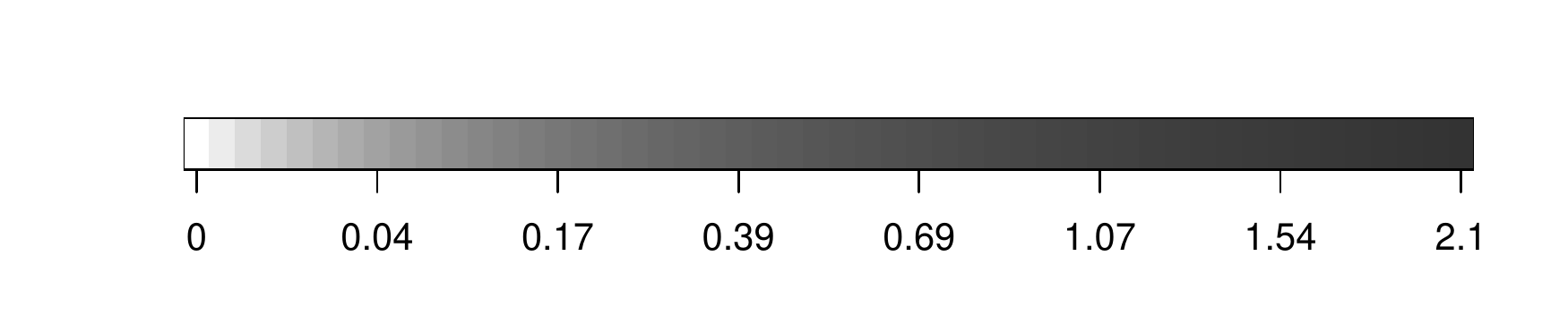}
	\caption{Plots of the observed neuron orientations (dots) and the fitted mixture density $\hat{f}$ (grey scale).  Left: plot of orientations represented as cosine-colatitude and longitude along with their marginal histograms and marginal fitted densities (solid curves) found by numerically integrating $\hat{f}$. Right: stereographic projection of the northern (top) and southern (bottom) hemisphere.}
	\label{fig:fitted_mixture_dist_vs_observed}
\end{figure}

Figure~\ref{fig:brain_cells} reveals no inhomogeneity for the neuron locations, whereas it is evident that the orientations are inhomogeneous pointing mostly toward the pial surface of the brain (the plane perpendicular to the $z$-axis). 
Thus, we estimated the intensity of the locations by
$\hat{\rho}_1 = n/|W| = 1.9 \times 10^{-5}$, and for the orientations we let the intensity be
$\rho_2(u)= c f(u)$, where $c$ is a positive parameter and $f$ is a density on $\mathbb{S}^2$ which
we model as follows.  Figure~\ref{fig:fitted_mixture_dist_vs_observed}
indicates that the orientations arise from a mixture of two
distributions; one distribution with points falling close to the North
Pole and another with points falling in a narrow girdle. Therefore, we
let $f(u) = pf_K(u) + (1 - p)f_W(u)$ be the mixture density of a Kent
and a Watson distribution on $\mathbb{S}^2$ \citep[see e.g.][for a
detailed description of these spherical
distributions]{Fisher:etal:1987}. In brief, the Kent density, $f_K$,
depends on five parameters (three directional, one concentration, and
one ovalness parameter), and its contours are oval with centre and
form specified by the directional parameters. Depending on the values
of the ovalness and concentration parameter, the Kent distribution is
either uni- or bimodal. Here, to account for the large number of
points centred around the North Pole, we consider the unimodal Kent
distribution. Furthermore, the Watson density, $f_W$, depends on two
parameters; a directional parameter determining the centres of the
density's circular contours, and a concentration parameter controlling
where and how fast the density peaks. Depending on the sign of the
concentration parameter, the density either decreases or increases as
the geodesic distance to the centres of the contours increases, giving
rise to either a bimodal or girdle shaped distribution. Since the
Watson distribution shall describe the orientations on the girdle, the
concentration parameter must be negative.

The nine parameters of the proposed intensity function $\rho_2$ were estimated as follows. We naturally estimate $c$ by $n$. The orientations occurring on the southern hemisphere are presumed to come from the Watson distribution, while the orientations on the northern hemisphere come from both distributions.  Therefore, and because the Watson density on the northern hemisphere is a  reflection of the southern hemisphere, we simply estimated the mixture probability by $\hat{p} = 1 - 2n_s/n = 0.94$, where $n_s$ is the observed number of points on the southern hemisphere. The directional parameters were chosen based on expectations expressed by the scientist behind the data collections, which were supported by visual inspection of the data; e.g.\ the directional parameter for the Kent distribution that determines the centre of the contours was chosen as the North Pole corresponding to the direction perpendicular to the pial surface and consistent with Figure~\ref{fig:fitted_mixture_dist_vs_observed}. Finally, the concentration and ovalness parameters were estimated by numerical maximization of their likelihood function, giving the estimated density
\begin{align*}
\hat{f}(u) = 0.94\, C_K \exp\{14.89 u_3 + 2.69(u_1^2 - u_2^2)\} +
0.06\, C_W\exp(-7.88 u_2^2),
\end{align*}
where $u = (u_1, u_2, u_3) \in \mathbb{S}^2$ and $C_K$, $C_W$ are normalising constants \citep[see][for details]{Fisher:etal:1987}. 
Figure~\ref{fig:fitted_mixture_dist_vs_observed} suggests that the fitted density (and associated marginal densities found by numerical integration of $\hat{f}$) adequately describe the distribution of the observed orientations. Therefore, we now turn to investigate whether the locations and orientations can be described by Poisson models with the estimated intensities. 

\begin{figure}[htpp]
	\centering
	\includegraphics[width=.4\textwidth]{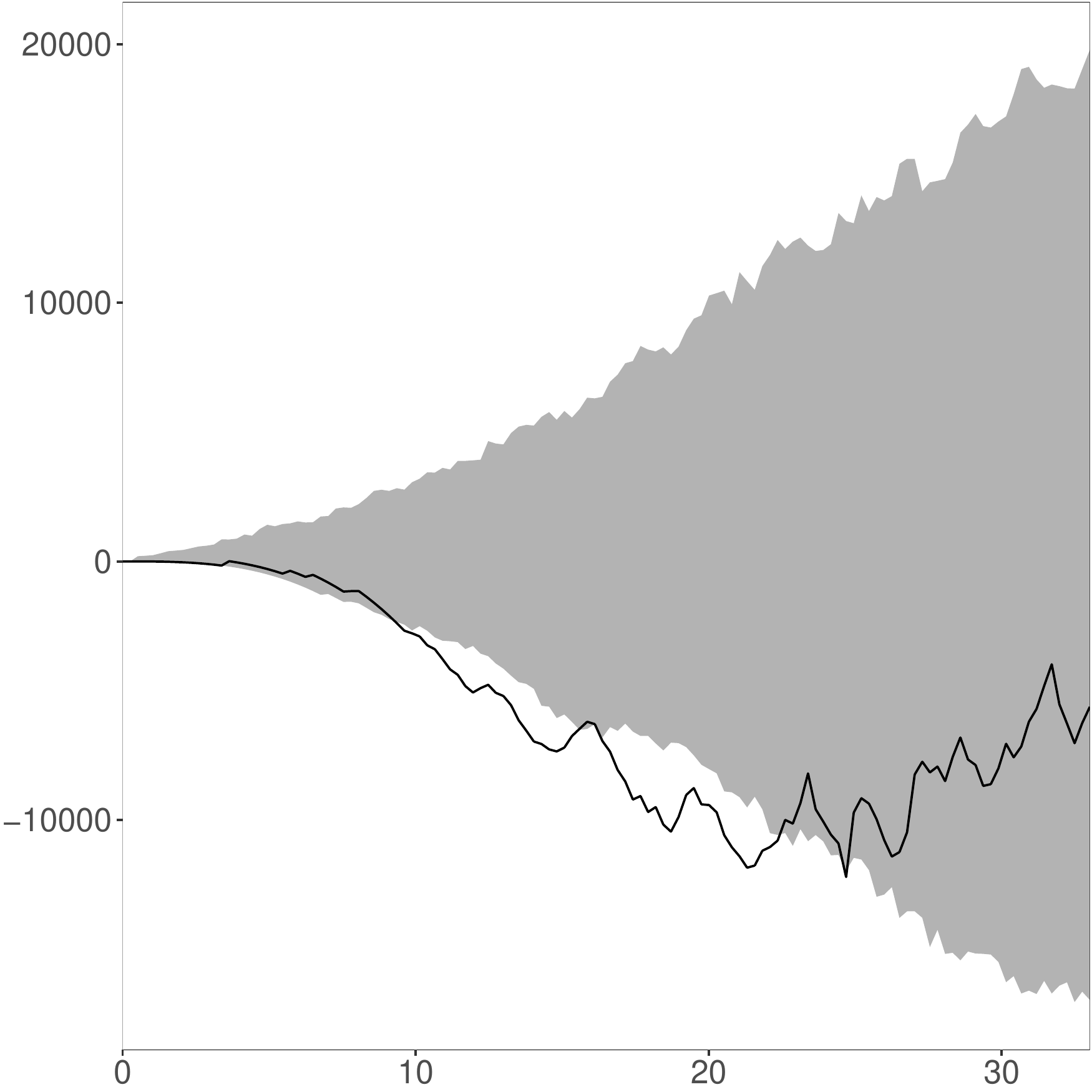}\hspace{5mm} 
	\includegraphics[width=.4\textwidth]{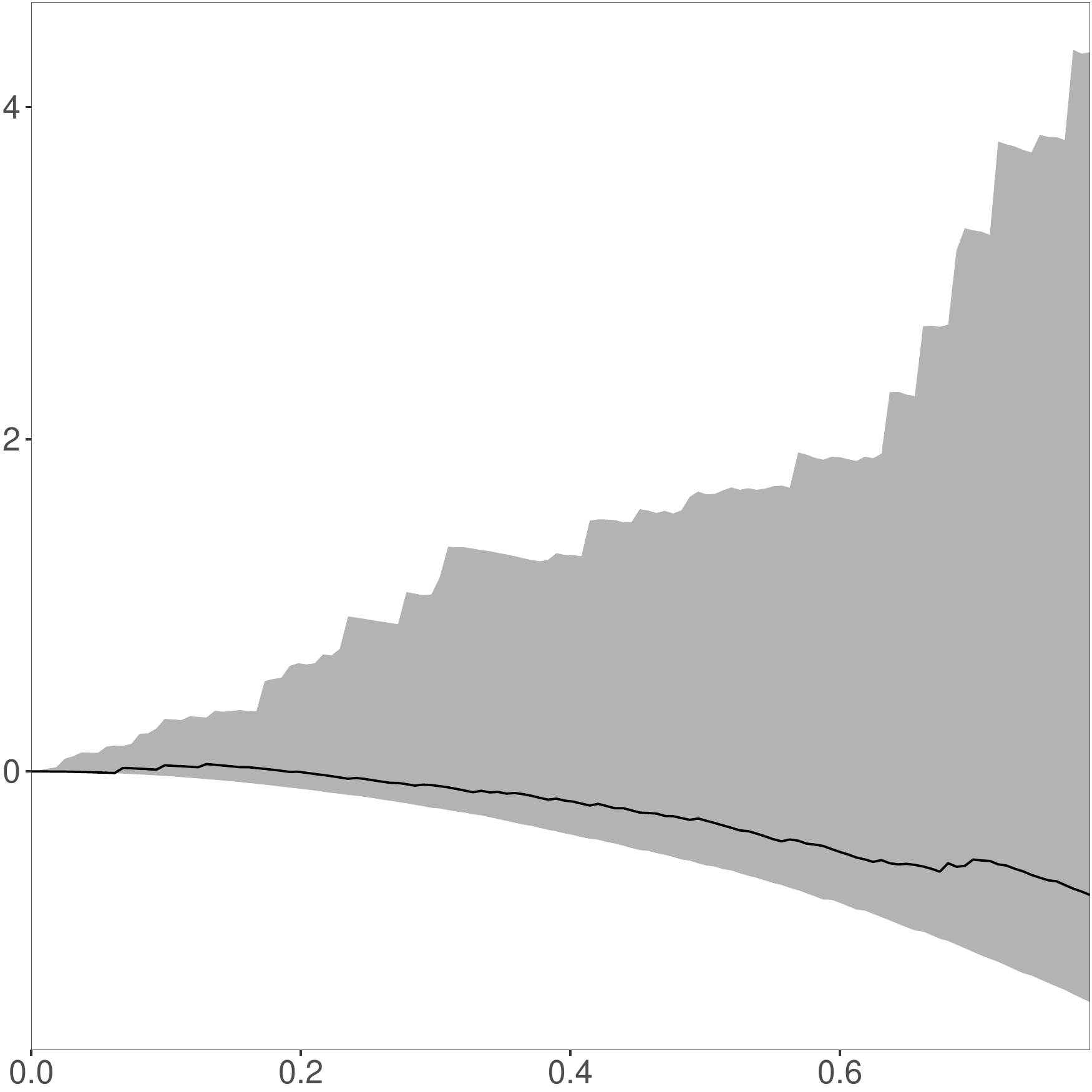}
	\caption{Left: $\hat{K}_{1}(r) - 4\pi r^3/3$ for the neuron locations (solid curve) along with a 95\% global rank envelope (grey area) under a homogeneous Poisson model on $W \subseteq \mathbb{R}^3$. Right: $\hat{K}_{2}(s) - 2\pi \{1 - \cos(s)\}$ for the neuron orientations (solid curve) along with a 95\% global rank envelope (grey area) under the fitted inhomogeneous Poisson model on $\mathbb{S}^{2}$.}
	\label{fig:marginal_K_under_Poisson}
\end{figure}

First, we considered the locations and orientations separately and
used $\hat{K}_1$ and~$\hat{K}_2$, respectively, as test functions for
the global rank envelope procedure. Using 2499 simulations from a
homogeneous Poisson process on $W$, we obtained a global rank envelope
test with $p$-interval $(0.000, 0.020)$. Thus, we reject that the
locations follow a homogeneous Poisson model. The associated $95\%$
global rank envelopes in Figure~\ref{fig:marginal_K_under_Poisson}
show that the rejection is due to the observed $\hat{K}_1(r)$ falling
below the envelope for $r$-values between
\SIrange{10}{25}{\micro\meter}. This suggests that the observed
locations exhibit some degree of regularity that needs to be
modelled. For the orientations, a global rank envelope test based on
2499 simulations under an inhomogeneous Poisson model on
$\mathbb{S}^2$ with intensity $\hat{\rho}_2(u) = n \hat{f}(u)$ gave a
$p$-interval of $(0.475, 0.481)$ and thus no evidence against the
proposed model. Figure~\ref{fig:marginal_K_under_Poisson} shows the
associated $95 \%$ global rank envelope.

\begin{figure}[!ht]
	\centering
	\includegraphics[width=.4\textwidth]{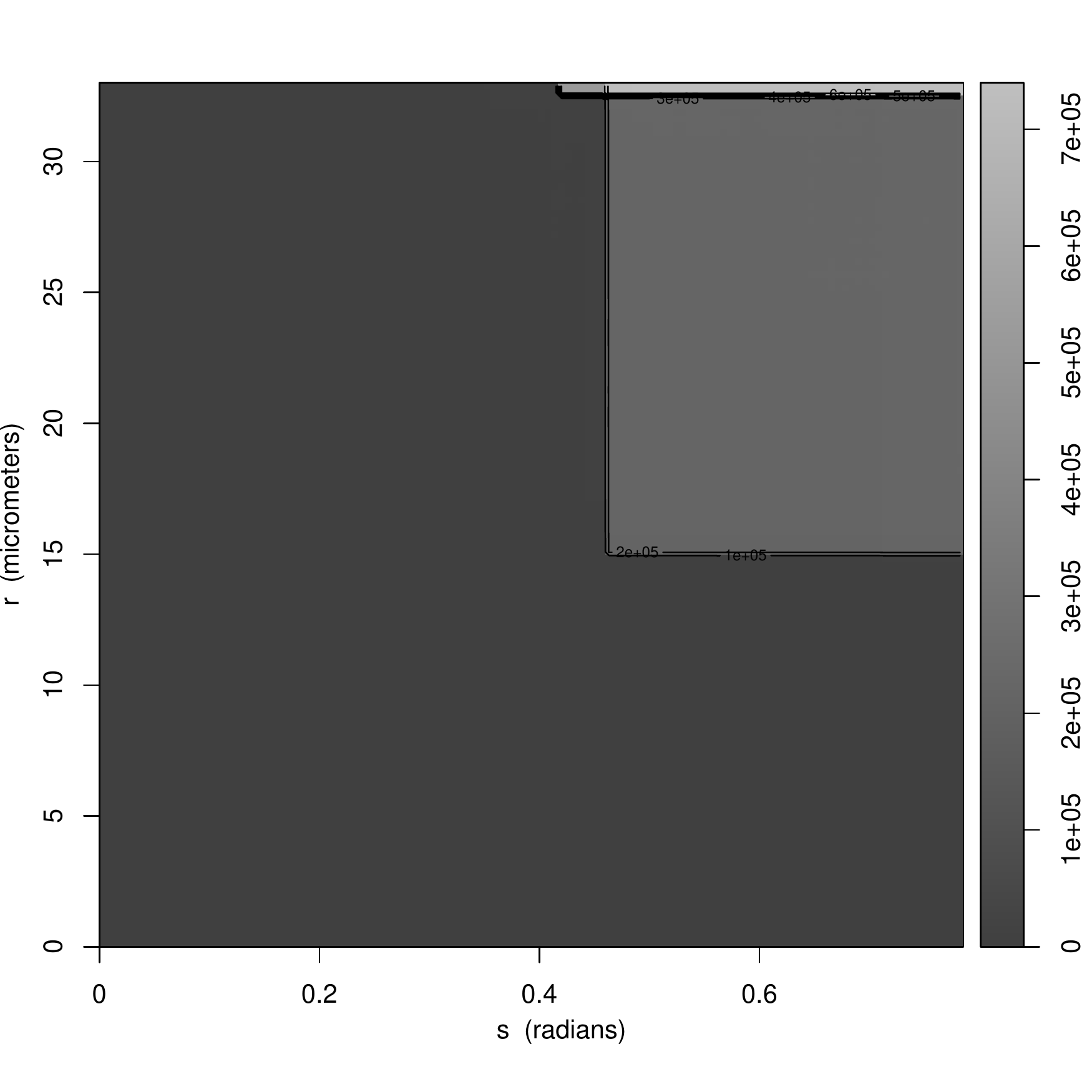}
	\includegraphics[width=.4\textwidth]{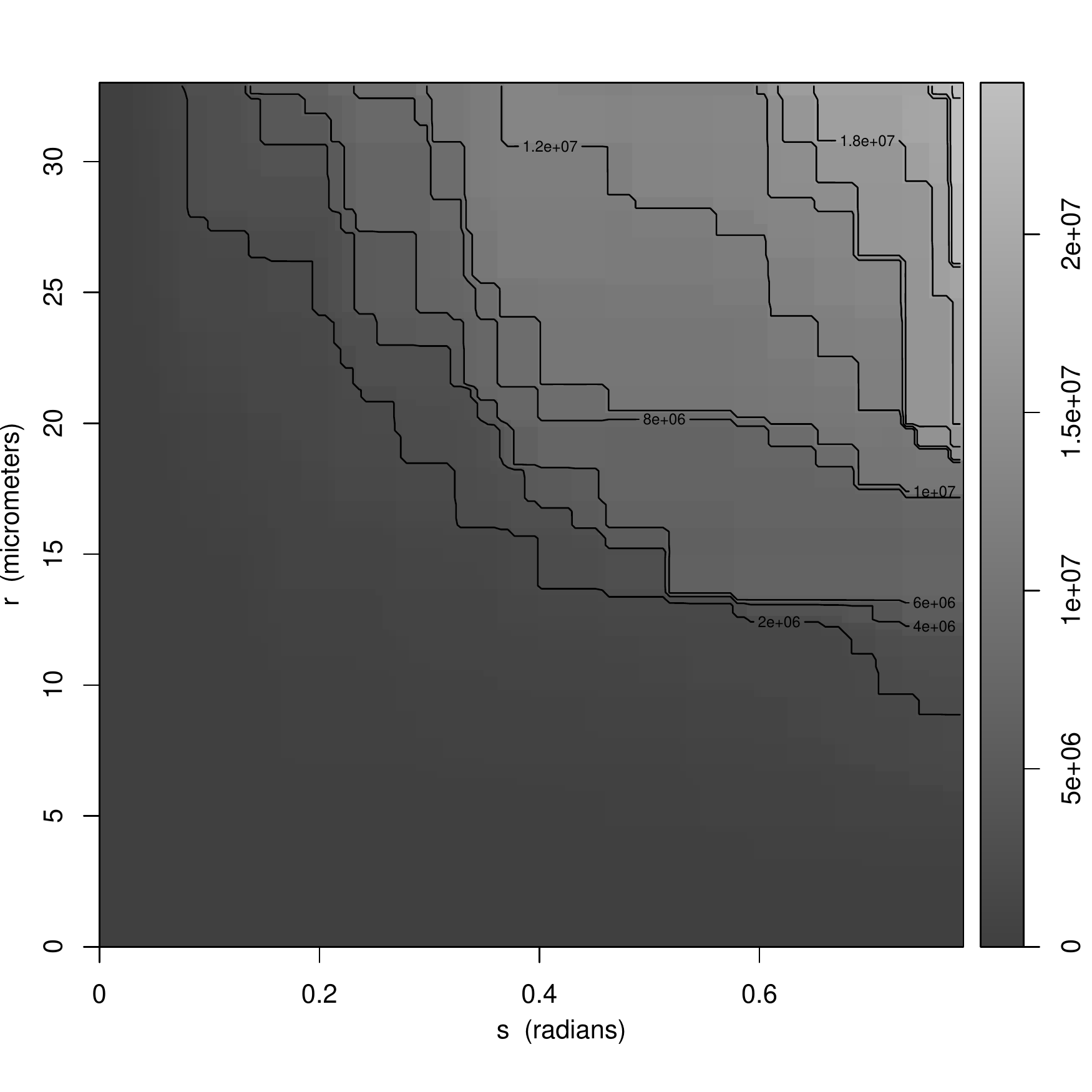}
	\\%
	\includegraphics[width=.4\textwidth]{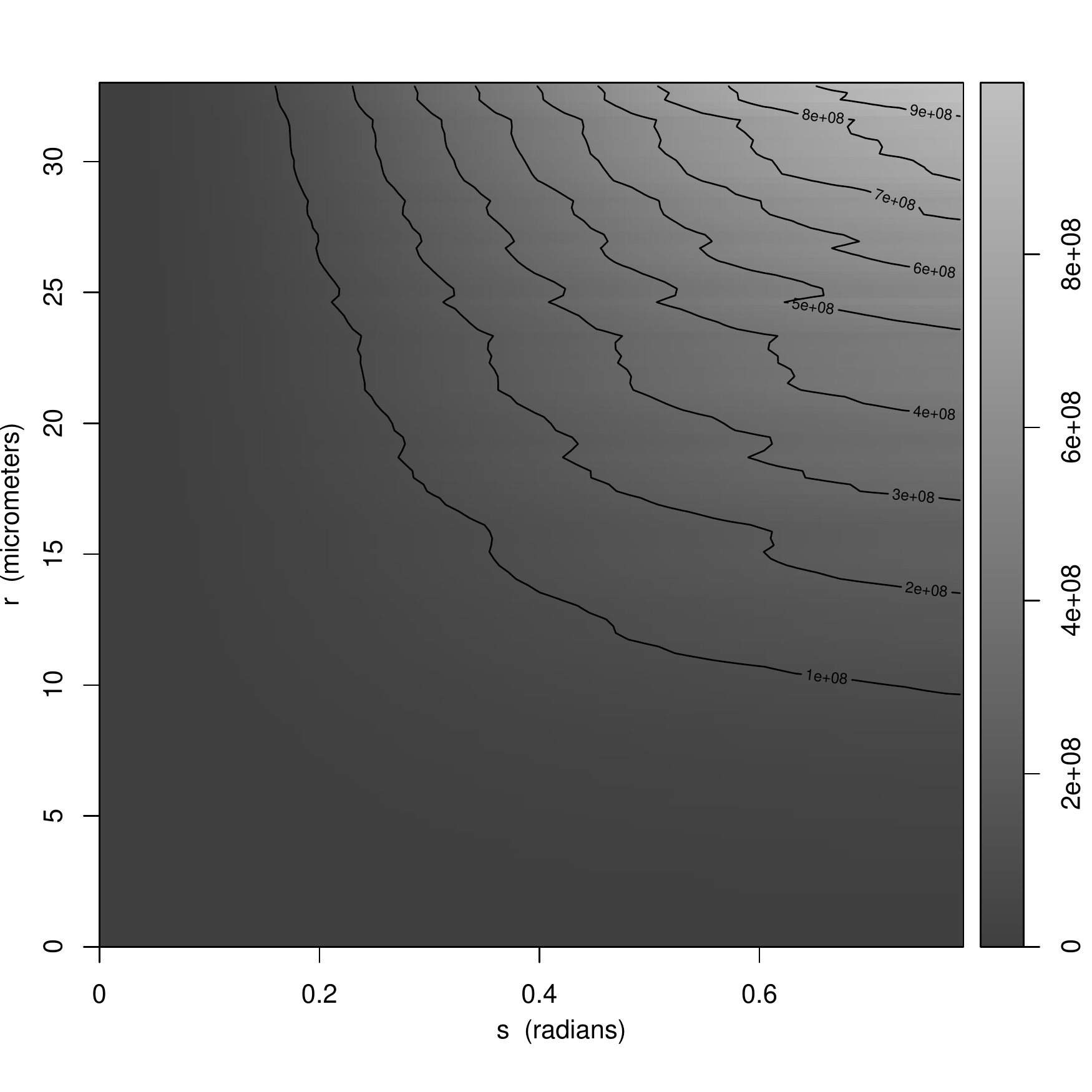}
	\includegraphics[width=.4\textwidth]{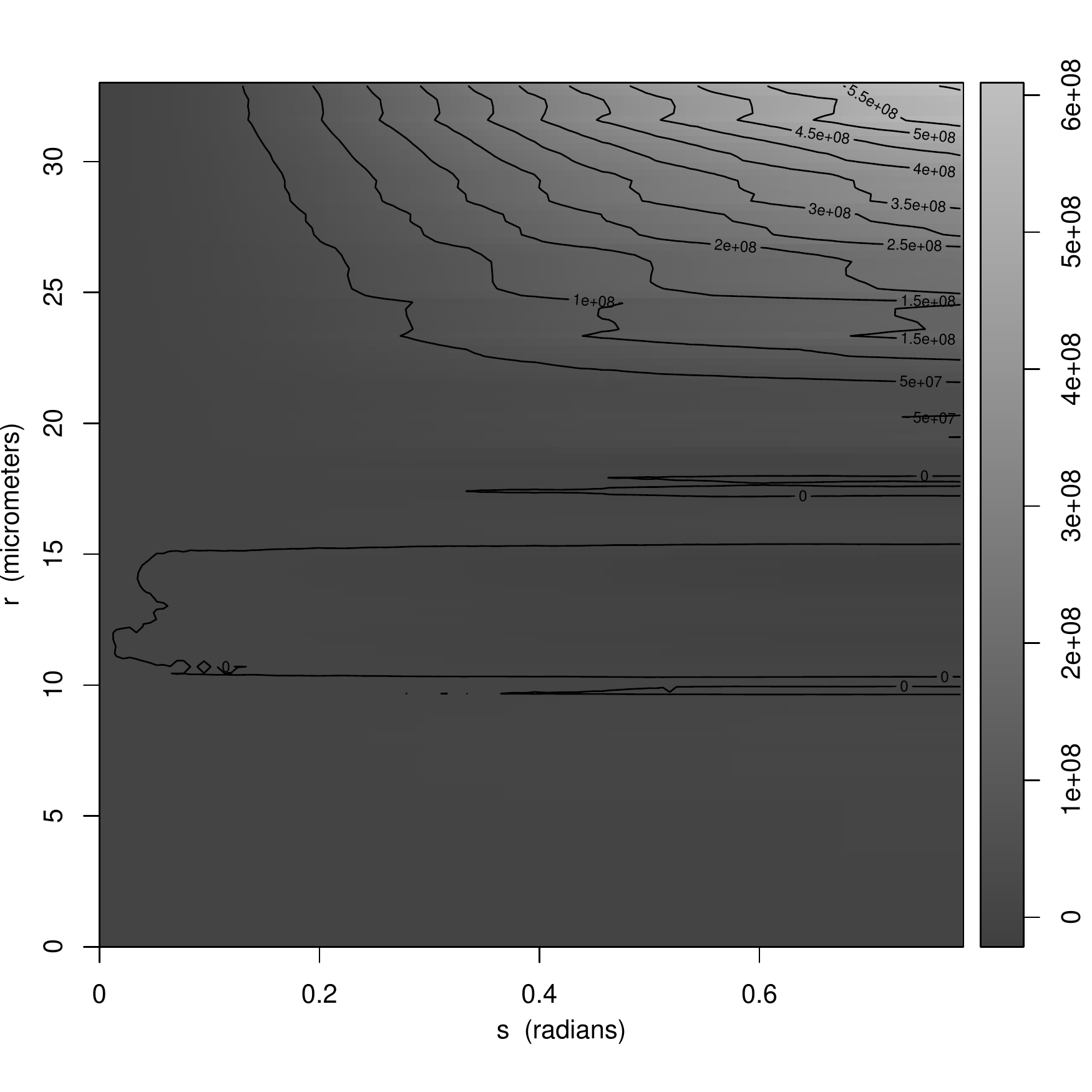}
	\caption{For $\hat{T} = \hat{K}$ (first row) and $\hat{T} = \hat{D}$ (second row): difference between $\hat{T}(r, s)$ for the observed neuron locations and orientations and the lower ($\hat{T}_\mathrm{low}$) or upper ($\hat{T}_\mathrm{upp}$) 95\% global rank envelope under the fitted inhomogeneous Poisson model on $W \times \mathbb{S}^2$. Left: $\hat{T}(r, s) - \hat{T}_\mathrm{low}(r, s)$.  Right: $\hat{T}_\mathrm{upp}(r, s) - \hat{T}(r, s)$.}
	\label{fig:KD_alldata_Poissontest}
\end{figure}

Second, we considered the data as a space-sphere point pattern and
used $\hat{K}$ and $\hat{D}$ to test for a Poisson model with a
separable intensity estimated by
$\hat{\rho}(y, u) = \hat{\rho}_1\hat{\rho}_2(u)/n = n\hat{f}(u)/|W|$,
cf.\ \eqref{e:sep1}. As the test functions $\hat{K}(r,s)$ and
$\hat{D}(r,s)$ depend on the two-dimensional argument $(r,s)$ and they
are non-smooth with large jumps due to few orientations occurring in
areas with low intensity, we increased the number of simulations to
49999 in order to improve the quality of the global rank envelope test
\citep[49999 is a much higher number than recommended in][for test
functions depending on one argument only]{Myllymaki2017}.  This
resulted in the $p$-intervals $(0.547, 0.549)$ for $\hat{K}$ and
$(0.000, 0.003)$ for $\hat{D}$; plots of the difference between the
associated envelopes and the observed test function are shown in
Figure~\ref{fig:KD_alldata_Poissontest}. In conclusion, the test based
on $\hat{K}$ reveals no evidence against the proposed space-sphere
Poisson model even though the corresponding Poisson model for the
locations was rejected by the test based on $\hat{K}_1$. However, the
test based on $\hat{D}$ provides a great deal of evidence against the
model. This conclusion is probably due to the fact that for this data
set $\hat{K}_1(r)\hat{K}_2(s) \gg \hat{K}(r, s)$, meaning that the
test based on $\hat{D}$ is highly controlled by the values of
$\hat{K}_1$ and $\hat{K}_2$, which cf.\ Figure~\ref{fig:KD_alldata_Poissontest} results in a rejection for
$r$-values from \SIrange{10}{20}{\micro\meter}, in line with the test
based on $\hat{K}_1$.

\goodbreak

It is unsatisfactory that $\hat{K}$ does not detect any deviation from Poisson when $\hat{K}_1$ clearly does, but we expect that the large jumps in $\hat{K}(r, s)$, caused by $(r, s)$-close neighbours with low intensity, may explain why no evidence against the model is detected. 
The few orientations that were modelled using a Watson distribution mostly fall in places with very low intensity. Therefore, we independently thinned the space-sphere point pattern with retention probability $\hat{p}\hat{f}_K(u) / \hat{f}(u)$ for $u \in \mathbb{S}^2$. Thereby (with high probability) we removed neurons with orientations that were most likely generated by the Watson distribution. This lead to removal of 26 neurons.
For the thinned data, the global rank envelope test based on $\hat{K}$ for testing the hypothesis of an inhomogeneous Poisson process with intensity proportional to a Kent density gave a $p$-interval of $(0.052, 0.058)$. Still, the model was not rejected at a $5\%$ significance level, but we at least got closer to a rejection; and so we continued the analysis with the thinned data. The analysis here indicates that, at least in some cases, the power of the global rank envelope test based on $\hat{K}$ may be small. This is investigated further in Section~\ref{s:simulationstudy}.

As we have seen, a homogeneous Poisson model is not adequate for the locations, and thus a Poisson model with intensity $\hat{\rho}$ as described above is not suitable for describing the space-sphere point pattern. To investigate whether orientations and locations can be modelled separately, that is, whether the locations and orientations are independent, we kept the locations fixed, and independent of the locations we simulated IID orientations from the fitted Kent distribution. The resulting global rank envelope test based on $49999$ of such simulations gave a $p$-interval of $(0.9255, 0.9258)$ for $\hat{K}$ and $(0.1265, 0.1266)$ for $\hat{D}$, showing no evidence against the hypothesis of independence between locations and orientations.  
Alternatively, if one does not have a suitable model to simulate the spherical (or spatial) components from, the independence test may be performed by randomly permuting the components. Formally, this tests only the hypothesis of exchangeability; a property that is fulfilled under independence. Performing such a permutation test for our data where the locations were fixed and the orientations permuted $49999$ times resulted in a $p$-interval of $(0.5431, 0.5445)$ using either $\hat{K}$ or $\hat{D}$ (as $\hat{K}$ and $\hat{D}$ only differ by a constant under permutation of the orientations and thus lead to equivalent tests).

\FloatBarrier

\subsection{Simulation study}\label{s:simulationstudy} 
In the data analyses in Sections~\ref{sec:fireball}--\ref{s:neuron}, the tests based on $\hat{K}$ failed to reject the proposed Poisson models in cases where the corresponding spatial model was rejected when using $\hat{K}_1$. 
To investigate whether the space-sphere $K$-function is a valuable
addition to the existing functional summary statistics on the space
and sphere, we performed a simulation study comparing the power of
global rank envelope tests based on either $\hat{K}$, $\hat{D}$, or a
combination of $\hat{K}_1$ and $\hat{K}_2$. (The combined test
function is simply a concatenation of $\hat{K}_1$ and $\hat{K}_2$. \cite{mrkvicka-etal-17} recommended using such a combination rather than $\hat{K}_1$ or $\hat{K}_2$ as a test function.) 

Specifically, we consider a homogeneous LGCP $X$ driven by a random
field $\Lambda(y, u) = \rho\exp\{Z(y,u)\}$, where $\rho > 0$ and
\begin{align*}
Z(y,u) = \alpha + \sigma_{1} Z_{1}(y) + \sigma_{2}Z_{2}(u) + \delta Z_{3}(y,u), \qquad y \in \mathbb{R},\ u \in \mathbb{S}^2,
\end{align*}
for parameters $\sigma_{1},\sigma_{2} > 0$, $\delta \geq 0$, and
$\alpha = -(\sigma_{1}^{2} + \sigma_{2}^{2} + \delta^{2})/2$. Further,
$Z_{1}, Z_{2}$, and $Z_{3}$ are independent GRFs with mean $0$ and
covariance functions
$c_{1}(y_{1},y_{2}) = \exp(- \|y_{1} -\nobreak y_{2}\|_d/\phi_{1})$,
$c_{2}(u_{1},u_{2}) = \exp(-d(u_{1},u_{2})/\phi_{2})$, and
$c_3(y_{1}, u_{1} ,y_{2}, u_{2}) = \linebreak
c_1(y_1, y_2)c_2(u_1,
u_2)$, respectively, with parameters $\phi_1, \phi_2 > 0$.  Note that
the resulting LGCP is homogeneous (and thus first order separable) and
SOIRS for any value of $\delta\ge0$. In addition, by \eqref{e:lgcp2},
the process is second order separable if and only if $\delta = 0$, in
which case $X$ has pair correlation function
\begin{equation}\label{eq:sim_pcf}
  \begin{aligned}
    \MoveEqLeft g_\theta(y_1, u_1, y_2, u_2)
    \\
    &= \exp\{ \sigma_1^2 c_1(y_1, y_2) + \sigma_2^2 c_2(u_1, u_2) \},
    \qquad y_1, y_2 \in \mathbb{R}, \ u_1, u_2 \in
    \mathbb{S}^2,
  \end{aligned}
\end{equation}
where $\theta = (\sigma_1, \phi_1, \sigma_2, \phi_2)$. 

For each value of $\delta = 0, 0.5, 1, 1.5, 2$, we simulated 100 realisations of a LGCP on $[0, 1] \times \mathbb{S}^2$ with $\rho = 1000$, $\sigma_1 = \sigma_{2} = 0.5$, $\phi_1 = 0.05$, and $\phi_2 = 0.132$.
Then for each of these simulations, we fitted the LGCP model with $\delta = 0$ using a second order composite likelihood approach proposed by \cite{Guan2006} to estimate $\theta$. 
In the present time-sphere setting, for a finite point pattern $x \subset [0, 1] \times \mathbb{S}^2$, the log second order  composite likelihood is given by
\begin{equation}\label{eq:CL}
  \begin{aligned}
    \mathrm{CL}(\theta; x)
    ={} &\sum_{(y_i, u_i), (y_j, u_j) \in x}^{\neq}
    w(y_i, u_i,y_j, u_j)
    \log\{\rho\psup{2}_\theta(y_i, u_i, y_j, u_j)\}
    \\
    &- n_{r, s}
    \begin{aligned}[t]
      \log \biggl\{& \int_{[0, 1] \times \mathbb{S}^2}\int_{[0,
        1]\times \mathbb{S}^2 }w(y_1, u_1,y_2,
      u_2)
      \\
      &\cdot\rho\psup{2}_\theta(y_1, u_1, y_2, u_2)\,\mathrm d\mu(y_1,
      u_1)\,\mathrm d\mu(y_2, u_2)\biggr\}.
    \end{aligned}
  \end{aligned}
\end{equation}
Here, for user specified distances $r$ and $s$,
$w(y_1, u_1,y_2, u_2) = \mathbb{I}\{\|y_1 - y_2\|_d < r$,
$d(u_1, u_2) <
s\}$, 
$n_{r, s}$ is the number of $(r, s)$-close neighbours, and
$\rho_\theta\psup{2}$ is the second order joint intensity function,
which for the homogeneous LGCP presented above is
$\rho_\theta\psup{2}(y_1, u_1, y_2, u_2) = \rho^2 g_\theta (y_1, u_1,
y_2, u_2)$. Then \eqref{eq:CL} is easily seen not to depend on $\rho$,
and by \eqref{eq:sim_pcf} the composite likelihood can be written as
\begin{align}\label{eq:CL_sep}
&\mathrm{CL}(\theta; x) = l_{1}(\sigma_1, \phi_1; x) + l_{2}(\sigma_2, \phi_2; x)
\end{align}
for functions $l_1$ and $l_2$. Thus, maximising the composite likelihood with respect to $\theta$ can be split into two maximisation problems; that is, maximising $l_1$ with respect to $(\sigma_1, \phi_1)$ and $l_2$ with respect to $(\sigma_2, \phi_2)$. 
Finally, we tested the null hypothesis $\delta = 0$ using the global rank envelope test with 4999 simulations from the fitted model using either $\hat{K}$, $\hat{D}$, or a combination of $\hat{K}_1$ and $\hat{K}_2$ as test functions. 

Table \ref{tab:simulation1} gives an overview of the conclusions reached by these tests. Note that the power of the tests based on either of the three test functions in general increases with $\delta$, both for the liberal and conservative test (for details see Appendix A). Thus, with increasing degree of non-separability the tests more often detect deviation from the separable model. However, tests based on $\hat{K}$ and particularly $\hat{D}$ seem preferable in this setup as they have a higher power than tests based on $\hat{K}_1$ combined with $\hat{K}_2$. 
 
Obviously, the conservative $p$-value always lead to fewer rejections than the liberal, giving a lower power. However, if the global rank envelope procedure is based on a higher number of simulations, then the conservative and liberal test will more often lead to the same conclusion.

\begin{table}
  \centering
  \newcommand\mc[1]{\multicolumn{1}{c}{#1}}
  \caption{Power of tests for different values of $\delta$ when using
    the global rank envelope test with either $\hat{K}$, $\hat{D}$, or
    $\hat{K}_1$ combined with $\hat{K}_2$. The decision was made using
    a significance level of $5 \%$ for both the liberal and
    conservative tests.} 
	\label{tab:simulation1}
	\begin{tabular}{
            ll
            S[table-format=2.0,table-space-text-post=\%]  <{\%}
            S[table-format=2.0,table-space-text-post=\%]  <{\%} 
            S[table-format=2.0,table-space-text-post=\%]  <{\%}
            S[table-format=2.0,table-space-text-post=\%]  <{\%}
            S[table-format=3.0,table-space-text-post=\%]  <{\%}
            l@{\kern-\tabcolsep}
          }
          \toprule
		& Test function & \mc{$\delta = 0$} & \mc{$\delta = 0.5$} & \mc{$\delta = 1$}   & \mc{$\delta = 1.5$}  & \mc{$\delta = 2$} \\
		\midrule
		Liberal & $\hat{K}$ & 4  & 7  & 42 & 75 & 98& \\
		& $\hat{D}$ & 2  & 45 & 92 & 97 & 100& \\
		& $\hat{K}_1, \hat{K}_2$ & 10 & 11 & 29 & 28 &
                42 & \\
                \midrule
		Conservative & $\hat{K}$ & 2  & 5  & 32 & 72 & 90& \\
		& $\hat{D}$ & 0  & 26 & 77 & 82 & 86 & \\
		& $\hat{K}_1, \hat{K}_2$ & 10 & 11 & 29 & 28 & 40 & \\
                \bottomrule
              \end{tabular}
\end{table}

\section{Additional comments}\label{s:discussion}
Section~\ref{sec:examples} introduced examples of space-sphere point processes for which the second order separability property described in Section~\ref{s:separability} seems natural. 
However, for other classes of point processes a different structure of the pair correlation function may be more interesting. For example, suppose $X$ is a Cox process driven by 
\begin{align}\label{e:SNCP_RF}
\Lambda(y,u) = \sum_{(y', u', \gamma ') \in \Phi}\gamma' k(y', u', y, u), \qquad y \in \mathbb{R}^d, \, u \in \mathbb{S}^k,
\end{align}
where $\Phi$ is a Poisson process on $S \times (0,\infty)$ with intensity function $\zeta$, and $k(y', u',\cdot,  \cdot)$ is a density with respect to $\mu$. 
Then $X$ is called a \textit{shot noise Cox process  (SNCP)} with kernel $k$ \citep{moeller:03}. The process has intensity function 
\[
\rho(y, u) = \iint \gamma' \zeta(y', u', \gamma')k(y', u', y, u)\, \mathrm d\mu(y', u') \, \mathrm d\gamma', \qquad y \in \mathbb R^d, \, u \in \mathbb{S}^k, 
\]
and pair correlation function
\begin{equation}\label{e:gSNCP}
  \begin{aligned}
    \MoveEqLeft g(y_1, u_1, y_2, u_2)
    \\
    &= 1 + \frac{\displaystyle \iint \gamma'^2 \zeta(y', u', \gamma')k(y', u', y_1,
      u_1)k(y', u', y_2, u_2) \, \mathrm d\mu(y', u') \, \mathrm
      d\gamma'}{\rho(y_1, u_1)\rho(y_2, u_2)}
  \end{aligned}
\end{equation}
for any $y_1, y_2 \in \mathbb{R}^d$ and any $u_1, u_2 \in \mathbb{S}^k$ with $\rho(y_1, u_1)\rho(y_2, u_2) > 0$. 
In the trivial case where the kernel $k(y', u', y, u)$ in \eqref{e:SNCP_RF} does not depend on $u$ (or $y$), the SNCP is both first and second order separable, with intensity and pair correlation functions that do not depend on the spherical (or spatial) components, and the process thus fulfils second order separability in the sense of \eqref{e:sep2}. 
However, the specific structure of the pair correlation function for a SNCP in \eqref{e:gSNCP} makes it more natural to look for a product structure in $g-1$ rather than $g$. That is, we may say that $X$ is second order separable if there exist Borel functions $h_1$ and $h_2$ such that 
\begin{align}\label{eq:2nd_separability_alternative}
g(y_1, u_1, y_2, u_2) - 1= h_1(y_1, y_2)h_2(u_1, u_2), \qquad y_1, y_2 \in \mathbb{R}^d, \, u_1, u_2 \in \mathbb{S}^k. 
\end{align}
This property is naturally fulfilled whenever we consider a Poisson
process or any marked point process with marks that are IID and
independent of the ground process as described in
Example~\ref{ex:IndMark}.  Now, think of $\Phi$ in \eqref{e:SNCP_RF}
as a marked point process with ground process
$\{(y, u) : (y,u,\gamma) \in \Phi\}$ and marks
$\{\gamma :(y,u,\gamma) \in \Phi \}$, and assume that the ground
process and the marks are independent processes, the ground process is
a homogeneous Poisson process on $S$ with intensity $\alpha > 0$, and
the marks are IID with mean $m_1$ and second moment $m_2$. If in
addition $k(y', u', y, u) = k_0\{y-y',\allowbreak d(u, u')\}$, then $\Lambda$ and
thus $X$ is stationary in space and isotropic on the sphere. Further,
$X$ is homogeneous with intensity $\rho= m_1\alpha$ and pair
correlation function
\begin{align}\label{e:gSNCP_hom}
g(y_1, u_1, y_2, u_2) = 1 + \frac{m_2}{\alpha m_1^2}\int k_0\{y_1-y', d(u_1, u')\}k_0\{y_2 - y', d(u_2, u')\} \, \mathrm d\mu(y',u')
\end{align}
for $ y_1, y_2 \in \mathbb R^d$ and $u_1, u_2 \in \mathbb S^k$. 
Clearly, \eqref{e:gSNCP_hom} depends only on $(y_1, y_2)$ through $y_1 - y_2$, and on $(u_1, u_2)$ through $d(u_1, u_2)$, although there is no simple expression for these dependencies in general.
Furthermore, separability in the form of \eqref{eq:2nd_separability_alternative} is fulfilled if the kernel $k$ in \eqref{e:gSNCP_hom} factorizes such that 
\begin{align*}
k_0\{y - y', d(u, u')\} = k_{01}(y-y')k_{02}\{d(u, u')\}, \qquad y, y' \in \mathbb R^d, \, u, u' \in \mathbb S^k,
\end{align*}
for Borel functions $k_{01}$ and $k_{02}$. Then by \eqref{e:rho1}--\eqref{e:g2}, the pair correlations functions for $Y$ and $U_W$ are
\begin{align*}
g_1(y_1, y_2) &= 1 + c_1\frac{m_2}{\alpha m_1^2} \int k_{01}(y_1-y')k_{01}(y_2-y')\,\mathrm dy', \qquad y_1, y_2 \in \mathbb R^d,
\end{align*}
and
\begin{align*}
g_2(u_1, u_2) &= 1 + c_2\frac{m_2}{\alpha m_1^2} \int k_{02}\{d(u_1, u')\}k_{02}\{d(u_2, u')\}\,\mathrm d\nu(u'), \qquad u_1, u_2 \in \mathbb S^k,
\end{align*}
where
\begin{align*}
  c_1 &= \frac{1}{\sigma_k^2} \iiint k_{02}\{d(u_1,
  u')\}k_{02}\{d(u_1, u')\}\,\mathrm d\nu(u_1)\,\mathrm
  d\nu(u_2)\,\mathrm d\nu(u')
\end{align*}
and
\begin{align*}
c_2 &= \frac{1}{|W|^2} \int \int_W \int_W k_{01}(y_1 - y')k_{01}(y_1 - y') \,\mathrm dy_1\,\mathrm dy_2\,\mathrm dy'.
\end{align*} 
That is, $g_1$ depends only on the spherical components through the
constant $c_1$, and similarly $g_2$ depends only on the spatial
components through $c_2$. \cite{DP2014} discuss how an analogue
property can be exploited to estimate the parameters of space-time
SNCPs using minimum contrast estimation for the projected processes. A
similar procedure will be applicable for space-sphere point processes,
but we have not investigated this further.

Section~\ref{s:est} considered the situation where the spatial
components of $X$ are observed on a subset of $\mathbb{R}^d$ and the
spherical components are observable on the entire sphere. In more
general applications, the spherical components may only be observable
on a subset of $\mathbb S^k$ leading to edge effects on the sphere
too. To account for this, edge correction methods for the sphere
should be used when estimating $K_2$ \mbox{\citep[see][]{LBMN2016}}
and $K$. If $X$ is observable on a product space $W_1 \times W_2$,
where $W_1 \subset \mathbb R^d$ and $W_2 \subset \mathbb S^k$, then an
edge corrected estimate for $K$ may be obtained by combining edge
corrected estimates for $K_1$ and $K_2$ analogous to
\eqref{e:estK}. Concerning the specific choice of edge correction
method, \cite{BRT2015} mentioned for planar point processes that, ``So
long as some kind of edge correction is performed \dots, the particular
choice of edge correction technique is usually not critical''. We
expect that the situation is similar for our setting.

For one-dimensional test functions, \cite{Myllymaki2017} recommend using 2499 simulations to perform a global rank envelope test, and \cite{mrkvicka-etal-17} discuss the appropriate number of simulations when using a multivariate test function (as the empirical space-sphere $K$-function). In Section~\ref{s:neuron}, we used 49999 simulations for the global rank envelope test based on $\hat{K}$, since $\hat{K}$ had steep jumps. 
To avoid this large number of simulations, a refinement of the global rank envelope test discussed in \cite{MHM2018} can be applied. 

In Example~\ref{ex:PoissonCox} we noticed that if a space-sphere point process is a Poisson process, then the spatial and spherical components are Poisson processes as well. Nevertheless,  using $\hat{K}$ for testing a space-sphere Poisson model  may lead to a different conclusion than using $\hat{K}_1$ and $\hat{K}_2$ for testing whether the corresponding Poisson models for the spatial and spherical components, respectively, are appropriate. Indeed, in the case of Figures~\ref{fig:MarginalKfun}--\ref{fig:KEnvelope}, the test based on $\hat{K}_1$ showed some evidence against a homogeneous Poisson model for the fireball event times, while no evidence against a homogeneous Poisson model for the locations over time was seen with the test based on $\hat{K}$. This observation together with the results in Section~\ref{s:neuron} was our motivation for making the simulation study in Section~\ref{s:simulationstudy}, where we investigated the power of global rank envelope tests based on either $\hat K$, $\hat D$, or a combination of $\hat K_1$ and $\hat K_2$, and we concluded that tests based on $\hat K$ and in particular $\hat D$ seem preferable.  


In Section~\ref{s:simulationstudy}, we utilised homogeneity and second order separability to speed up optimisation of the second order composite likelihood proposed in \cite{Guan2006}. If the process is inhomogeneous but first (and second) order separable, we still get a separation of the composite likelihood similar to $\eqref{eq:CL_sep}$, where $l_1$ and $l_2$ now may depend also on intensity parameters.
As an alternative, the second order composite likelihood discussed in \cite{Waagepetersen2007} can be used. However, in that case, first and second order separability do not yield a separable likelihood as in \eqref{eq:CL_sep}, and for our simulation study it resulted in unstable estimates \citep[and thus it was discarded in favour of the one proposed by][]{Guan2006}. Furthermore, one may investigate whether the adaptive procedure discussed in \cite{Lavancier2018} will provide stable estimates in the space-sphere setting. In short, \cite{Lavancier2018} consider the score function related to \eqref{eq:CL} and introduce a modified weight function $w$ depending on $g$.

In this paper, we considered point processes living on $\mathbb{R}^d\times \mathbb{S}^k$. Naturally, we may extend the results/methods to more general metric spaces $\mathbb{R}^d \times \mathbb{M}$, where $\mathbb{M}$ is a compact set (e.g.\ a torus). However, we need to require some invariance property for the metric space $\mathbb{M}$ and its metric under a group action, such that we can define an equivalence of the SOIRS property needed to define $K$.

\paragraph{Acknowledgements.}
The authors are grateful to Ji\v{r}{\'i} Dvo\v{r}{\'a}k for helpful
comments and to Ali H.\ Rafati for collecting the pyramidal cell
data. This work was supported by The Danish Council for Independent
Research | Natural Sciences, grant DFF -- 7014-00074 ``Statistics for
point processes in space and beyond'', and by the ``Centre for
Stochastic Geometry and Advanced Bioimaging'', funded by grant 8721
from the Villum Foundation.

\section*{Appendix A}

In Sections~\ref{sec:fireball}--\ref{s:simulationstudy}, we used the global rank envelope test presented in \cite{Myllymaki2017} to test for various point process models. In this appendix, we briefly explain the idea and use of such a test. A global rank envelope test compares a chosen test function for the observed data with the distribution of the test function under the null model; as this distribution is typically unknown it is approximated using a Monte Carlo approach.
The comparison is based on a rank that only gives a weak ordering of the test functions. Thus, instead of a single $p$-value, the global rank envelope test provides
an interval of $p$-values, where the end points specify the most liberal and conservative $p$-values of the test. 
A narrow $p$-interval is desirable as the test is inconclusive if the $p$-interval contains the chosen significance level. 
The width of the $p$-interval depends on the number of simulations, smoothness of the test functions and  dimensionality. \cite{Myllymaki2017} recommended to use 2499 simulations for one-dimensional test functions and a significance level of $5\%$. 

An advantage of the global rank envelope procedure is that it provides a graphical interpretation of the test in form of critical bounds (called a global rank envelope) for the test function. For example, if the observed test function is not completely inside the $95 \%$ global rank envelope, this corresponds to a rejection of the null hypothesis at a significance level of $5\%$. Furthermore, locations where the observed test function falls outside the global rank envelope reveal possible reasons for rejecting the null model. 

In their supplementary material, \cite{Myllymaki2017} discussed two approaches for calculating test functions that rely on an estimate of the intensity.  One approach is to reuse the intensity estimate for the observed point pattern in calculation of all the test functions, another is to reestimate the intensity for each simulation and then use this estimate when calculating the associated test function.  
For the $L$-function, which is a transformation of $K_1$, \cite{Myllymaki2017} concluded that the reestimation approach give the more powerful test. In this paper, we have therefore based all our global rank envelope tests on that approach. 

\bibliographystyle{spbasic}      
\bibliography{references}   

\end{document}